\documentclass[11pt]{amsart}
\usepackage[english]{babel}
\usepackage{a4wide,color,graphicx,amsmath,amssymb,verbatim,mathrsfs,latexsym,bm,dsfont}
\usepackage{mathtools}
\allowdisplaybreaks
\newtheorem{theo}{Theorem}[section]
\newtheorem{prop}[theo]{Proposition}
\newtheorem{lemma}[theo]{Lemma}

\newtheorem{rem}[theo]{Remark}

\newtheorem{defn}[theo]{Definition}
\newcommand{\R}{\mathbb{R}}	%
\newcommand{\N}{\mathbb{N}}	%
\newcommand{\pa}{\partial}		%
\newcommand{\Div}{\textrm{div}\,}	%

\newcommand{\fun}[5]{  			%
	\begin{array}{cccc}
		{#1}\,: & #2 & \to & #3 \\
		\phantom{{#1}\,:\,} & #4 & \mapsto & #5
\end{array} }

\newcommand{\chf}[1]{{\raisebox{3pt}{\Large $\chi$}}_{#1}}
\newcommand{\vel}{\bm{v}}

\newcommand{\rhob}{\overline{\rho}}
\newcommand{\dx}{\,\mbox{d}x}
\newcommand{\dt}{\,\mbox{d}t}
\newcommand{\bJ}{\bm{J}}
\newcommand{\bM}{\mathbb{M}}
\newcommand{\bI}{\mathbb{I}}
\newcommand{\bS}{\mathbb{S}}
\newcommand{\bT}{\mathbb{T}}
\newcommand{\bW}{\mathbb{W}}
\newcommand{\dd}{\,\mbox{d}}
\newcommand{\dom}{\mathscr{D}}
\newcommand{\nfa}{\textrm{NFA}}
\newcommand{\bx}{\bm{x}}

\title[Modeling and analysis of a phase-field model]{Modeling and analysis of a phase-field model for organic solar cells}
\author{Jan-Frederik Pietschmann}
\address{Universit\"{a}t Augsburg, Institut f\"ur Mathematik and  Centre for Advanced Analytics and Predictive Sciences (CAAPS), Universit\"{a}tsstra\ss e 12a, 86159 Augsburg, Germany}
\email{jan-f.pietschmann@uni-a.de}
\author{Carmen Tretmans}
\address{Universit\"{a}t Augsburg, Institut f\"ur Mathematik, Universit\"{a}tsstra\ss e 12a, 86159 Augsburg, Germany}
\email{carmen.tretmans@uni-a.de}
\author{Nicola Zamponi}
\address{Universit\"{a}t Augsburg, Institut f\"ur Mathematik, Universit\"{a}tsstra\ss e 12a, 86159 Augsburg, Germany}
\email{nicola.zamponi@uni-a.de}

\date{\today}

\begin{document}
	
    \begin{abstract}
    We derive and analyze a phase-field model describing the formation of the active layer of organic solar cells. The process combines phase separation of an electron donor (polymer) and an electron acceptor (non-fullerene acceptor) with evaporation of a solvent. Taking the macroscopic motion of the fluid into account, this results in a coupled system of incompressible Navier--Stokes, degenerate multispecies Cahn--Hilliard, and Allen--Cahn equations with cross-diffusion terms. We give a thermodynamically consistent derivation of the system and prove the existence of global-in-time weak solutions without any size constraints on the data. The solutions dissipate the total energy and require an ad-hoc variational formulation of the definition of the diffusion fluxes, which is made necessary by the degeneracy of the mobility matrix. 
	\end{abstract}

	\maketitle
	
	\tableofcontents

    \textbf{Use of AI:} AI was used for language checking and consistency of notation. After completing the first draft of this paper, we used AI (Opus~5.5, ChatGPT~Astra, accessed Nov.~2026) to review the mathematical arguments. This raised several technical issues, which we resolved. The final results were checked by all authors and we take full responsibility for the content of this work.  \newline

    \textbf{Acknowledgment: } The authors are grateful to Olivier J.\ J.\ Ronsin and his group for insightful discussions and for generously sharing their expertise on organic photovoltaic devices. NZ thanks the University of Augsburg for financial support via the research program “Forschungspotentiale besser nutzen!”. JFP and CT thank the Deutsche Forschungsgemeinschaft (DFG) for funding this work (Research Unit FOR 5387 POPULAR, project no.\ 461909888).
    
	\section{Introduction}\label{sec.intro}
    Organic photovoltaic cells (OPVs) have emerged as a promising low-cost, sustainable, and flexible alternative to inorganic (e.g., silicon-based) cells~\cite{Zhang2021}. While there are still shortcomings in terms of efficiency and stability~\cite{Jiang2024, Liu2026}, the possibility of cost-effective mass production by printing techniques makes them an active area of research~\cite{Zhang2022, Dong2025}. OPVs are typically fabricated by solvent-based thin-film deposition, using methods ranging from spin coating for lab-scale devices to roll-to-roll methods for large-scale production~\cite{Taranenko2026}. In each case, an electron donor (polymer) and an electron acceptor (non-fullerene acceptor) are dissolved in at least one solvent to form a homogeneous mixture~\cite{Hou2018}. The mixture is then coated onto a substrate, creating a thin film. As the solvent evaporates, the donor and acceptor undergo phase separation, a process that halts once the solvent is depleted. This results in a film consisting of donor- and acceptor-rich regions, forming the active layer of the solar cell. The donor--acceptor morphology is highly sensitive to processing conditions, such as blend composition, evaporation rate, and substrate characteristics~\cite{Ciloglu25}. In turn, the morphology largely determines device performance~\cite{Huang2014}. A clear understanding of how processing conditions govern morphology formation is therefore essential for achieving high device efficiency. 
    
    To accurately model the morphology formation, several coupled physical phenomena must be accounted for: (i) phase separation of the donor and acceptor materials, (ii) evaporation of the solvent, and (iii) the hydrodynamic effects of the blend. These effects are closely coupled. For instance, the amount of solvent remaining in the system directly governs the rate of phase separation. 

    Various models addressing different aspects of these phenomena have been proposed. The phase separation of donor and acceptor materials within organic solar cells is described, e.g., in~\cite{WoGa12, Bergermann23} by means of phase-field methods. In these works, evaporation is governed by an outflow boundary condition for the solvent. Ronsin et al.~\cite{RJEBH21} and Ronsin and Harting~\cite{RonHar22} include evaporation effects by means of an additional order parameter, driven by an energy functional. Describing phase separation including hydrodynamic effects requires the formulation of coupled Cahn--Hilliard Navier--Stokes equations. These systems date back to \textit{Model~H} by Hohenberg and Halperin~\cite{ModelH}, followed by the model of Lowengrub and Truskinovsky~\cite{LowengrubTruskinovsky}, which allows for non-matching densities. Later, Abels et al.~\cite{AbGaGr2012} presented a model introducing the divergence-free volume-averaged velocity. It is worth noting the work by ten Eikelder et al.~\cite{Eikelder2023} and ten Eikelder~\cite{Eikelder2025}, which presents a comprehensive overview of existing Cahn--Hilliard Navier--Stokes models. 
    
    In this work, we present a model describing all three effects: phase separation, solvent evaporation, and hydrodynamic effects. For the first two effects, we adopt the methods proposed in Ronsin et al.~\cite{RJEBH21} and Ronsin and Harting~\cite{RonHar22}, whereas for the latter we apply the framework of Abels et al.~\cite{AbGaGr2012}. More precisely, our model 
    \begin{itemize}
        \item describes donor--acceptor demixing by a phase-field method. This leads to Cahn--Hilliard-type equations, driven by a free energy,
        \item treats solvent evaporation through an order parameter evolving according to an Allen--Cahn equation. Likewise, evaporation is driven by a free energy,
        \item accounts for hydrodynamic effects through the Navier--Stokes equations,
        \item is thermodynamically consistent, i.e., it satisfies both the first and second laws of thermodynamics.
    \end{itemize}
    The resulting model takes the form of coupled Cahn--Hilliard, Navier--Stokes, and Allen--Cahn equations.
    
    Analytical problems related to coupled Cahn--Hilliard and Navier--Stokes equations have already been investigated, in particular in the case of single-component fluids without cross-diffusion effects. In~\cite{AbeDepGar13_diff_dens}, the existence of weak solutions to a model developed in~\cite{AbGaGr2012} is proved. The model employs a diffuse-interface phase-field approach to describe the flow of two viscous incompressible Newtonian fluids, with the phase field evolving according to a Cahn--Hilliard equation. The resulting coupled Cahn--Hilliard Navier--Stokes system is thermodynamically consistent and dissipates the total energy as a consequence. The analysis is extended in~\cite{AbeDepGar13_degen_mob} to the case of a degenerate mobility. The local-in-time existence of strong solutions for regular energy densities is shown in~\cite{AbeWeb21,Web16}, while the case of less regular energy densities is treated in~\cite{Gio21,Gio22}. The long-time convergence towards a steady state and the propagation of regularity are proved in~\cite{AbeGarGio24}. A Hele-Shaw Cahn--Hilliard model is studied in~\cite{DedGarLam18}, which is derived via formal asymptotics from the Navier--Stokes Cahn--Hilliard system presented in~\cite{AbGaGr2012}.
    In~\cite{GarKovTra22}, a phase-field model for tumor growth including viscoelastic effects is presented and studied, which consists of coupled Cahn--Hilliard equations, reaction-diffusion equation, Navier--Stokes equations, and a mechanical force balance. 
    In~\cite{EbeGarNue20}, a phase-field model for tumor growth is introduced, which consists of a coupled Cahn--Hilliard equation for the tumor tissue,  a
    reaction-diffusion-advection equation for the nutrients, and a Brinkman--Stokes-type law for the fluid velocity. The authors prove global-in-time existence of weak solutions for the model with degenerate mobility.
    
    Cahn--Hilliard equations for multicomponent systems with cross-diffusion have also been studied.
    In~\cite{EhrMarPie21}, a Cahn--Hilliard model for a multicomponent mixture with cross-diffusion effects and degenerate mobility is studied. The authors prove global-in-time existence of weak solutions by exploiting the formal gradient-flow structure of the system and a suitable extension of the boundedness-by-entropy method~\cite{BurgerDiFrancescoPietschmannetal.2010,Jue15}; the notion of weak solution is adapted to possible degeneracies. 
    Stationary solutions and large-time asymptotics for the model are studied in~\cite{CVEMP24}.
    In~\cite{GioHeWu26}, the authors prove global-in-time existence of weak solutions to a Navier--Stokes Cahn--Hilliard system for viscous incompressible two-phase flows that incorporates the effects of chemotaxis, mass transport, and reactions. The model
    presents a cross-diffusion structure similar to the Keller--Segel model for chemotaxis.
    In~\cite{JueLi24}, a model describing the pre-patterning of lymphatic vessel morphology in collagen gels is presented, which consists of degenerate Cahn--Hilliard cross-diffusion equations with a singular potential. The authors prove the global-in-time existence of weak solutions to the system by employing a three-level approximation scheme as well as energy and entropy inequalities.
    In~\cite{HopKinMueWag26}, a degenerate Cahn--Hilliard model for viscoelastic phase separation involving cross-diffusion coupling to a bulk stress variable is studied, with a focus on a formal sharp-interface limit.

    In summary, considerable work has been done on the analysis of either Navier--Stokes Cahn--Hilliard models for single-species systems with no cross effects or Cahn--Hilliard equations for multicomponent systems with cross-diffusion, whereas coupled Cahn--Hilliard Navier--Stokes systems for fluid mixtures with cross-diffusion have not received much attention in the analytical literature. In particular, to the best of our knowledge, the Cahn--Hilliard Navier--Stokes system presented and studied in this paper has neither been previously derived nor investigated in the existing literature.

    This paper is organized as follows. In Section~\ref{sec.model.deriv}, the governing equations describing the morphology formation within organic solar cells are derived in a thermodynamically consistent way. Subsequently, in Section~\ref{sec.analysis}, an existence result for is stated and proven, using a time-discretization together with suitable a-priori bounds. The existence of iterates of the approximate scheme is shown in Appendix~\ref{appendix.proof.prop}.

    \section{Model derivation}\label{sec.model.deriv}
    \subsection{Notation and assumptions} \label{ss.assumptions}
    Let $\Omega\subset\R^3$ be a bounded domain with boundary $\partial\Omega\in C^{1,1}$ and outer unit normal $\nu$. The domain represents the drying active layer consisting of a polymer acting as the electron donor, a non-fullerene acceptor (NFA) acting as the electron acceptor, and a solvent, as well as an ambient air layer above the active layer. For $T>0$, define $Q_T=\Omega\times (0,T)$.   
    We consider the evolution of polymer, NFA, solvent, and air, labeled by $i \in \mathcal{S}^*$, $\mathcal{S}^* = \{p, \nfa, s, a\}$. Each constituent is described by a phase field $0 \leq \phi_i \leq 1$, $i \in \mathcal{S}^*$, indicating its volume fraction. The volume fractions are subject to the \emph{volume constraint}
    \begin{align*}
        \sum_{i \in \mathcal{S}^*} \phi_i(\bx, t) = 1 \qquad \mbox{a.e.~in }Q_T. 
    \end{align*} 
    This relation can be used to eliminate one constituent, and we additionally define the subset $\mathcal{S} = \{p, \nfa, s\}$, noting that $\phi_a = 1-\sum_{i \in \mathcal{S}} \phi_i$. Furthermore, we introduce the notations
    \begin{align*}
        \phi_{\mathcal{S}} = (\phi_p, \phi_\nfa,\phi_s) \text{ and } \phi_\mathcal{S}^* = (\phi_p, \phi_\nfa,\phi_s, \phi_a),  
    \end{align*} 
    and, analogously, $ \nabla \phi_{\mathcal{S}} = (\nabla \phi_p, \nabla \phi_\nfa,\nabla \phi_s) $ and $ \nabla \phi_\mathcal{S}^* = (\nabla \phi_p, \nabla \phi_\nfa,\nabla \phi_s, \nabla \phi_a)$.
    Finally, we denote the symmetric and skew-symmetric parts of the vector gradient by 
    \begin{align*}
		D(\vel) = \frac{1}{2}(\nabla\vel + \nabla\vel^\top) \text{ and } W(\vel) = \frac{1}{2}(\nabla\vel - \nabla\vel^\top),
	\end{align*}
    respectively. In what follows, we make various modeling choices and physical simplifications. Each constituent is assumed to be incompressible with constant density, regardless of whether it is in liquid or vapor state. Consequently, changes in density upon evaporation are neglected. We acknowledge that this is a significant simplification. However, since our focus is on the phase-separating behavior in the liquid phase, we believe this choice to be justified. Furthermore, the temperature $\theta$ is considered constant. Heat fluxes and heat sources (e.g., evaporative cooling) are neglected.
    
    \subsection{Conservation of mass}
    \label{sec.mass}
    All constituents are assumed to be incompressible with constant mass density $\rhob_i > 0, \ i \in \mathcal{S}^*,$ in the unmixed state, and volumes are assumed to be additive upon mixing. Consequently, the mass densities $\rho_i, \ i \in \mathcal{S}^*,$ are directly related to the constituents' volume fractions via $\phi_i(\bx, t) = \frac{\rho_i (\bx, t)}{\rhob_i}$. Each constituent moves with its own velocity $\vel_i(\bx, t)$ and satisfies an individual mass balance equation 
    \begin{align*}
      \partial_t \rho_i + \Div (\rho_i \vel_i) &= 0 \qquad \mbox{in }Q_T, \ i \in \mathcal{S}^*,  
    \end{align*}
    with mass densities $\rho_i(\bx, t)$, $i \in \mathcal{S}^*$. Equivalently,  
    \begin{align*}
        \partial_t \phi_i + \Div (\phi_i \vel_i) &= 0 \qquad \mbox{in }Q_T, \ i \in \mathcal{S}^*. 
    \end{align*}
    The mass density of the mixture is
    \begin{align*}
        \rho(\bx, t) &= \sum_{i \in \mathcal{S}^*}\rho_i(\bx, t) = \sum_{i \in \mathcal{S}^*} \rhob_i \phi_i(\bx, t) 
        = \rhob_a + \sum_{i \in \mathcal{S}} (\rhob_i - \rhob_a)\phi_i(\bx, t) . 
    \end{align*}
    To describe the mixture as a whole, one may introduce either a mass-averaged velocity or a volume-averaged velocity field, see the derivations in~\cite{AbGaGr2012, Eikelder2023}. We adopt the latter to ensure a divergence-free velocity field and define the volume-averaged velocity
    \begin{align*}
        \vel(\bx, t) = \sum_{i \in \mathcal{S}^*} \phi_i(\bx, t) \vel_i (\bx, t) .  
    \end{align*}
    Indeed, it follows that
    \begin{align}
        \Div \vel & = \Div \left( \sum_{i\in \mathcal{S}^*} \phi_i \vel_i \right) =  - \partial_t \left(\sum_{ i\in \mathcal{S}^* }  \phi_i \right) = 0  \qquad \mbox{in }Q_T.
        \label{eq_incompressibility}
    \end{align}
    The motion of each individual constituent can be decomposed into the volume-averaged motion and the motion relative to this average. The individual mass balance equations can be rewritten as 
    \begin{align*}
        \partial_t \rho_i & = -\Div  (\rho_i \vel_i) = -\Div \left(   \rho_i (\vel + \vel_i - \vel) \right) =  - \Div \left( \rho_i \vel + \bm{j}_i \right) \qquad \mbox{in }Q_T, \ i \in \mathcal{S}^* ,
    \end{align*}
    where $\bm{j}_i = \rho_i \vel_i - \rho_i \vel$ denotes the mass flux of constituent $i$ relative to the volume-averaged velocity. Using $\Div \vel = 0 $ and rescaling by $\rhob_i$, we arrive at the formulation 
    \begin{align}
        \label{eq_volumebalance}
        \partial_t \phi_i + \vel \cdot \nabla \phi_i = - \Div \bm{J}_i \qquad \mbox{in }Q_T, \ i \in \mathcal{S}^*  , 
    \end{align}
    defining the volume fluxes $\bm{J}_i = \frac{1}{\rhob_i} \bm{j}_i$. By the volume constraint, $\sum_{i \in \mathcal{S}^*} \bm{J}_i = 0$. Summing the mass balance equations over all constituents yields the mass balance of the mixture,
    \begin{align}
        \partial_t \rho + \vel \cdot \nabla \rho + \Div \bm{j} = 0 \qquad \mbox{in }Q_T ,
        \label{eq_massbalance}
    \end{align}
    with total mass flux $\bm{j} = \sum_{i\in \mathcal{S}^*} \bm{j}_i$.
    
    \subsection{Conservation of momentum}
    \label{sec.momentum}
    We consider the balance of momentum of the macroscopic motion of the fluid. Following Gurtin et al.~\cite{GuPoVi1996}, we assume that the momenta of the constituents relative to the gross motion of the fluid are negligible. This allows us to treat the mixture as a single fluid with velocity $\vel$, whose momentum is approximated by $\rho \vel$ rather than by $\sum_{i \in \mathcal{S}^*} \rho_i \vel_i = \rho \vel + \sum_{i \in \mathcal{S}^*} \bm{j}_i$. Note that $\sum_{i \in \mathcal{S}^*} \bm{j}_i = \sum_{i \in \mathcal{S}}\frac{\rhob_i - \rhob_a}{\rhob_i} \rho_i \left(\vel_i - \vel \right)$ by the volume constraint. Hence, the approximation of the momentum is exact if all densities coincide. Under this assumption, following Cauchy's theorem, the local linear momentum balance reads 
    \begin{align*}
        \partial_t ( \rho \vel ) + \Div (\rho \vel \otimes \vel) = \Div (\tilde \bT) \qquad \mbox{in }Q_T
    \end{align*}
    for a surface force tensor $\tilde \bT$. Substituting the mass balance equation \eqref{eq_massbalance} and the incompressibility constraint \eqref{eq_incompressibility} yields 
    \begin{align}
    \label{eq_mombalance}
        \rho \partial_t \vel + (\nabla \vel )(\rho \vel + \bm{j}) = \Div(\tilde \bT + \vel \otimes \bm{j}) \qquad \mbox{in }Q_T. 
    \end{align}
    Following the derivation by Abels et al.~\cite{AbGaGr2012}, we define the symmetric objective stress tensor $\bT = \tilde \bT + \vel \otimes \bm{j}$, which is split into a viscous stress tensor $\bS$ and a volumetric part $- p \bI$, hence $\bT  = \bS - p \bI$. With this substitution, the local momentum balance reads 
    \begin{align}
       \rho \partial_t \vel + (\nabla \vel)(\rho \vel + \bm{j})  = \Div \bS - \nabla p \qquad \mbox{in }Q_T .
       \label{eq_momentumbalance}
    \end{align}
    For future reference, we state the local momentum balance in its divergence form, 
    \begin{align}\label{eq.mom.2}
		&\pa_t(\rho\vel) + \Div (\bW) + \nabla p = 0\quad\mbox{in }Q_T,
	\end{align}
    with $\bW = \rho\vel\otimes\vel + \vel\otimes \bm{j} - \bS$. 
    
    \subsection{Evaporation}
    \label{sec.evap}
    Solvent evaporation is governed by two coupled processes. First, solvent changes its state of matter from liquid to gaseous at the liquid--gas interface. Second, the gaseous solvent diffuses throughout the surrounding air. To describe the first phenomenon, we allow all constituents to be either in a liquid or vapor state (of matter). Following Ronsin et al.~\cite{RJEBH21} and Ronsin and Harting~\cite{RonHar22}, we introduce an order parameter $\Phi_v \colon Q_T\rightarrow \R$, with
    \begin{align*}
        0 \leq \Phi_v \leq 1 \qquad\mbox{a.e.~in }Q_T,
    \end{align*}
    denoting the volume fraction of the mixture in the vapor phase. Here, $\Phi_v$ describes the fraction of all constituents in the vapor phase simultaneously, and $(1-\Phi_v)$ the corresponding fraction of the mixture in the liquid phase. While we allow all constituents to be in either liquid or vapor phase, the liquid--vapor transition of the individual constituents will be controlled by their respective free energy contributions. As explained in Subsection~\ref{ss.assumptions}, although the constituents are allowed to have different densities $\rhob_i$, we assume these densities to be constant and independent of the constituents' state of matter. Although the order parameter $\Phi_v$ has no influence on the mixture's density, it does influence various other quantities, e.g., the system's energy and the molecular mobilities. The second phenomenon is described by diffusion of the gaseous solvent through the air, which is driven by thermodynamic forces. 
    \begin{rem} \label{rem_bc}
    In a real experiment, diffusion of solvent and air is not confined to $\Omega$, but continues across the domain boundary $\partial \Omega$ into the ambient air. In our setting, we assume that the liquid film only occupies a small part of the total domain, resulting in a large gaseous layer on top of the domain. More precisely, we assume sufficient air in our domain to absorb all solvent without saturating. However, if the domain size does not satisfy this requirement, due to, e.g., numerical restrictions, diffusion into the ambient air has to be modeled explicitly by an outflow boundary condition for solvent and a corresponding inflow boundary condition for air. Following Ronsin et al.~\cite{RJEBH21}, one can impose the flux boundary conditions on the (open part of the) boundary, $\Gamma_{\mathrm{open}} \subset \partial \Omega$, 
        \begin{align*}
            \bm{J}_s \cdot \nu = J_{\mathrm{out}} \  \text{ and } \ \bm{J}_a \cdot \nu = - J_{\mathrm{out}} \quad \text{on} \quad \Gamma_{\mathrm{open}} \times (0,T),  
        \end{align*}
        where $J_\mathrm{out}$ describes the outward volume flux of solvent and depends on the solvent volume fraction and the vapor order parameter. We refrain from any outflow boundary conditions and impose no-flux conditions, compare Subsection~\ref{ss.bc_ic}, 
        \begin{align*}
            \bm{J}_s \cdot \nu = 0 \  \text{ and } \ \bm{J}_a \cdot \nu = 0 \quad \text{on} \quad \partial \Omega \times (0,T).  
        \end{align*}    
    \end{rem}  
   
	\subsection{Thermodynamic consistency and constitutive equations} \label{sec.thdyn}
    It remains to define constitutive equations for the fluxes $\bm{j}_i$, the deviatoric stress tensor $\bS$, and an evolution equation for the vapor order parameter $\Phi_v$ compatible with the first and second laws of thermodynamics. Following the first law of thermodynamics, the total energy in an arbitrary control volume is conserved, i.e.,
    \begin{align*}
        \frac{\dd}{\dt} \int_{V(t)} e \dx = - \int_{\partial V(t)} \bm{j}_e \cdot \nu \dd s_x ,
    \end{align*}
    for the energy density per unit volume $e$ and energy flux $\bm{j}_e$. The energy is the sum of kinetic energy and internal energy. Approximating the kinetic energy by $\rho \frac{|\vel|^2}{2}$, and thereby neglecting the motion relative to the volume-averaged velocity, the energy density reads
    \begin{align*}
        e = \rho \frac{|\vel|^2}{2}  + u ,  
    \end{align*}
    where $u$ denotes the internal energy density per volume. Introducing the (Helmholtz) free energy density $f=  u - \theta s$ with entropy density per volume $s$ results in
    \begin{align*}
         \frac{\dd}{\dt} \int_{V(t)} \rho \frac{|\vel|^2}{2} +  f + \theta s \dx = - \int_{\partial V(t)} \bm{j}_e \cdot \nu \dd s_x .
    \end{align*}
    The second law of thermodynamics consists of the entropy balance and the entropy law 
    \begin{align*}
        \frac{\dd}{\dt} \int_{V(t)} \theta s \dx & = - \int_{\partial V(t)} \bm{j}_{\mathrm{ent}} \cdot \nu \dd s_x + \int_{V(t)} \sigma \dx , \\
        \nonumber
        \sigma & \geq 0 ,
    \end{align*}
    where $\bm{j}_{\mathrm{ent}}$ denotes $\theta$ times the entropy flux and $\sigma$ denotes $\theta$ times the entropy source per unit volume. Combining the first and second laws results in
    \begin{align}
    \label{modeling.dissip.ineq}
         \frac{\dd}{\dt} \int_{V(t)} \rho \frac{|\vel|^2}{2} + f   \dx =  - \int_{\partial V(t)} \bm{j}_e \cdot \nu \dd s_x + \int_{\partial V(t)} \bm{j}_{\mathrm{ent}} \cdot \nu \dd s_x  - \int_{V(t)} \sigma \dx ,
    \end{align}  
    with $\sigma \geq 0$. Combining the entropy law with equation~\eqref{modeling.dissip.ineq} results in the well-known \textit{dissipation inequality}. We opt to preserve the equation in equality form and will refer to it as the \textit{energy balance}.

    Various derivations have been presented to obtain suitable constitutive equations satisfying the energy balance and entropy constraint, see Abels et al.~\cite{AbGaGr2012}. We utilize the derivation based on a microforce balance, as proposed by Gurtin~\cite{Gurtin1996} in the form used by Abels et al.~\cite{AbGaGr2012}. We assume the free energy density $f$ to depend solely on the phase fields and their gradients, i.e., 
    \begin{align*}
        f(\phi_{\mathcal{S}^*}, \Phi_v, \nabla \phi_{\mathcal{S}^*}, \nabla  \Phi_v).    
    \end{align*}  
    Specifying the entropy and energy fluxes $\bm{j}_{\mathrm{ent}}$ and $\bm{j}_e$, and following the derivation presented in Appendix~\ref{app.thdyncons}, we find the following constitutive equations to be compatible with both the energy balance and the entropy constraint: 
    \begin{align}
    \label{eq_const_S}
        \bS & = 2 \eta(\phi_{\mathcal{S}^*})D (\vel) - \sum_{i \in \mathcal{S}^*} \nabla \phi_i \otimes \partial_{\nabla \phi_i} f  - \nabla \Phi_v \otimes \partial_{\nabla \Phi_v} f ,   \\
        \label{eq_const_j}
        \bm{j}_i & = - \sum_{j \in \mathcal{S}^*} m_{ij}(\phi_{\mathcal{S}^*}, \Phi_v) \nabla \mu_j ,  \\ 
         \label{eq_const_mu} 
        \mu_i &= \partial_{\phi_i} f  - \Div \partial_{\nabla \phi_i} f ,\\
        \label{eq_const_Phiv}
       	\pa_t\Phi_v + \vel\cdot\nabla\Phi_v	&= -m_v(\phi_{\mathcal{S}^*}) \left (\pa_{\Phi_v}f - \Div\,\pa_{\nabla\Phi_v}f \right). 
    \end{align}
    Here, we introduce the chemical potentials $\mu_i, i \in \mathcal{S}^*$. Furthermore, we define the viscosity $\eta(\phi_{\mathcal{S}^*}) \geq 0$, the symmetric positive semidefinite molecular mobility matrix $(\bM)_{ij} = \frac{m_{ij}(\phi_{\mathcal{S}^*}, \Phi_v)}{\rhob_i}$ satisfying $\sum_i \frac{m_{ij}}{\rhob_i} = 0 $, and the vapor mobility $m_v(\phi_{\mathcal{S}^*}) \geq 0$. Moreover, we require the free energy functional to obey $\partial_{\nabla \phi_i} f = \alpha_i \nabla \phi_i$ and $\partial_{\nabla \Phi_v} f = \alpha_v \nabla \Phi_v$ for some scalar (phase- and order-field-dependent) $\alpha_i$, $\alpha_v$, and lastly, the volume constraint to hold. The entropy production then reads
    \begin{align*}
		\sigma = 2\eta\,|D(\vel)|^2
		+ \sum_{i,j \in \mathcal{S}^*} \frac{m_{ij}}{\rhob_i} \,\nabla\mu_i\cdot\nabla\mu_j
		+ m_v \bigl|\partial_{\Phi_v} f - \Div \partial_{\nabla \Phi_v} f\bigr|^2
		\;\ge\;0 .
	\end{align*}
    The first and last terms of the entropy production are clearly nonnegative. Moreover, the positive semi-definiteness of the mobility matrix ensures that the middle term is nonnegative as well. Consequently, the constitutive equations satisfy the entropy law.

\subsection{Free energy density} \label{sec.freeenergy}
    The free energy of the system is given by $ F = \int_{\Omega} (f + f_0)  \dx $, where $f_0$ denotes a reference free energy density and $f$ the corresponding change in free energy density. Here, we derive the change in free energy density with respect to the unmixed liquid state of the constituents. The free energy density consists of a local free energy part $f^{\mathrm{loc}}$, describing the free energy in the bulk, and a nonlocal part $f^{\mathrm{nonloc}}$, accounting for spatial inhomogeneities. To find an expression for the local part, we first consider the liquid and vapor phases separately. In the liquid phase, the local free energy density consists of an entropic contribution due to mixing and an enthalpic contribution due to molecular interactions. Following the Flory--Huggins theory~\cite{Flory1942, Huggins1941}, we obtain 
    \begin{align*}
        f^{\mathrm{loc}}_{\mathrm{liq}}(\phi_{\mathcal{S}^*}) & = \frac{R \theta}{V^0}\left( \sum_{i \in \mathcal{S}^*} \frac{\phi_i}{N_i} \log \phi_i +  \sum_{\{i,j\} \subset \mathcal{S}^*}
        \chi_{ij} \phi_i \phi_j \right), 
    \end{align*}
     where $R$ denotes the gas constant, $V^{0}$ a reference molar volume, $N_i$ the degree of polymerization of constituent $i$, and $\chi_{ij}$ the Flory--Huggins interaction parameter between constituents $i$ and $j$, capturing the molecular interactions. In the vapor phase, the local free energy density consists of a contribution due to mixing, as well as a contribution due to the evaporation process. We assume molecular interactions to be negligible, hence
    \begin{align*}
        f^{\mathrm{loc}}_{\mathrm{vap}}(\phi_{\mathcal{S}^*}) & = \frac{R \theta}{V^0} \sum_{i \in \mathcal{S}^*} \left( \frac{\phi_i}{N_i} \log \phi_i - \frac{\phi_i}{N_i} \log \phi^{\mathrm{sat}}_i \right). 
    \end{align*}
    The first term describes the entropy of mixing, which, for $N_i = 1$, coincides with that of an ideal gas mixture. The second term describes the change in free energy due to evaporation, where $0 <\phi^{\mathrm{sat}}_i<1$ denotes the constant saturated vapor fraction. For the non-volatile polymer and NFA, ${\phi^{\mathrm{sat}}_p}$ and $ \phi^{\mathrm{sat}}_\nfa$ are chosen close to zero, while $\phi^\mathrm{sat}_a$ is chosen close to one, resulting in different liquid--vapor energy contributions. The value of $\phi^\mathrm{sat}_s$ is chosen intermediate and will largely govern the evaporation speed. The total change in the local free energy density is obtained by interpolating between liquid and vapor phases and is given by
    \begin{align*}
        f^{\mathrm{loc}}(\phi_{\mathcal{S}^*}, \Phi_v) &= \left( 1- \omega(\Phi_v)\right) f^{\mathrm{loc}}_{\mathrm{liq}}(\phi_{\mathcal{S}^*}) + \omega(\Phi_v) f^{\mathrm{loc}}_{\mathrm{vap}} (\phi_{\mathcal{S}^*}) . 
    \end{align*}
    Here, we introduce the interpolation function $\omega\in C^2(\R)$ satisfying $\omega=0$ on $(-\infty,0]$, $\omega=1$ on $[1,\infty)$, and $\omega$ increasing on $(0,1)$. The interfacial nonlocal term represents the contribution of surface tension and penalizes the field gradients, see the derivation in~\cite{CahnHilliard1958}, 
    \begin{align*}
        f^{\mathrm{nonloc}}(\nabla \phi_{\mathcal{S}^*}, \nabla \Phi_v) & = \sum_{i \in \mathcal{S}^*} \frac{\epsilon_i^2}{2} |\nabla \phi_i|^2 + \frac{\lambda^2}{2} |\nabla \Phi_v|^2,
    \end{align*}
    where $\epsilon_i$ and $\lambda$ denote surface tension parameters. Note that this expression fulfills $\partial_{\nabla \phi_i} f = \alpha_i \nabla \phi_i = \epsilon_i^2 \nabla \phi_i$. Altogether, the free energy density reads
    \begin{align}
    \label{eq_freeenergy}
        f(\phi_{\mathcal{S}^*}, \Phi_v, \nabla \phi_{\mathcal{S}^*}, \nabla \Phi_v) &= \frac{R \theta}{V^0} \left( 1- \omega(\Phi_v)\right) \left( \sum_{i \in \mathcal{S}^*} \frac{\phi_i}{N_i} \log \phi_i + \sum_{\{i,j\} \subset \mathcal{S}^*} \chi_{ij} \phi_i \phi_j  \right) \\ 
        & \quad +  \frac{R \theta}{V^0}  \omega(\Phi_v) \left(  \sum_{i \in \mathcal{S}^*} \frac{\phi_i}{N_i} \log \frac{\phi_i}{\phi^{\mathrm{sat}}_i}\right) \nonumber \\ 
        & \quad + \sum_{i \in \mathcal{S}^*} \frac{\epsilon_i^2}{2} |\nabla \phi_i|^2 + \frac{\lambda^2}{2}|\nabla \Phi_v|^2 . \nonumber
    \end{align}

    \subsection{Boundary and initial conditions} \label{ss.bc_ic}
    We impose no-flux boundary conditions for all constituents, compare Remark~\ref{rem_bc},
    \begin{align*}
        \left(\sum_{j \in \mathcal{S}^*} \frac{m_{ij}}{\rhob_i}(\phi_{\mathcal{S}^*}, \Phi_v) \nabla \mu_j\right) \cdot \nu  &= 0 \quad \text{on} \quad \partial \Omega \times (0,T), \ i \in \mathcal{S}^*.  
    \end{align*}
    The no-flux conditions are complemented with Neumann boundary conditions for the phase and order parameters
    \begin{align*}
        \partial_\nu \phi_i &= 0 \quad \text{on} \quad \partial \Omega \times (0,T),\ i \in \mathcal{S}^*, \\
        \partial_\nu \Phi_v &= 0 \quad \text{on} \quad \partial \Omega \times (0,T), 
    \end{align*}
    and the velocity field satisfies the free-slip condition
    \begin{align*}
        \vel \cdot \nu &= 0 \quad \text{on} \quad \partial \Omega \times (0,T), \\ 
        (\bI-\nu\otimes\nu)D(\vel)\nu &= 0  \quad \text{on} \quad \partial \Omega \times (0,T). 
    \end{align*}
    Furthermore, we prescribe the initial conditions
    \begin{align*}
		\vel(\cdot,0)=\vel^{\mathrm{in}},\quad\phi_i(\cdot,0)=\phi_{i}^{\mathrm{in}},\quad
		\Phi_v(\cdot,0) = \Phi_{v}^{\mathrm{in}}
		\qquad\mbox{in }\Omega, \quad i \in \mathcal{S}^*. 
	\end{align*}

    \subsection{Nondimensional governing equations}\label{sec.nondimeq}
    Here, we summarize the nondimensional form of the governing equations presented in the previous subsections. The full derivation can be found in Appendix~\ref{app.nondimeq}. The volume fractions are required to be nonnegative and sum to one, while the vapor order parameter $\Phi_v$ must be bounded between $0$ and $1$:
	\begin{align*}
		\phi_i \geq 0, \quad \sum_{i \in \mathcal{S}^*} \phi_i = 1 \quad \text{ and } \quad
        0\leq\Phi_v\leq 1\qquad\mbox{a.e.~in }Q_T.
	\end{align*}
    After choosing suitable length, time, and mass scales, the governing equations read
    \begin{align}\label{eq.v}
		&\rho\pa_t\vel + \left(\rho\vel 
		+ \sum_{j \in \mathcal{S}^*}\rhob_j \bJ_j\right)\cdot\nabla\vel \\ 
		\nonumber
		& =\Div\left(2\eta(\phi_{\mathcal{S}^*})D(\vel)\right)
		-\sum_{j \in \mathcal{S}^*} \beta_j\Div(\nabla\phi_j\otimes\nabla\phi_j)\\
		&-\beta_v\Div(\nabla\Phi_v\otimes\nabla\Phi_v) - \nabla p 
		& \mbox{in }Q_T, \notag\\
		\label{eq.divv}
		&\Div\vel = 0 & \mbox{in }Q_T,\\
		\label{eq.phii}
		&\pa_t\phi_i + \vel\cdot\nabla\phi_i + \Div \bJ_i = 0 ,
		\quad i \in \mathcal{S}^*,
		& \mbox{in }Q_T,\\
		\label{eq.Ji}
		&\bJ_i = -\sum_{j \in \mathcal{S}^*}\gamma_{ij}(\phi_{\mathcal{S}^*},\Phi_v)\nabla\mu_j,
		\quad i \in \mathcal{S}^*,
		& \mbox{in }Q_T,\\
		\label{eq.mui}
		&\mu_i = \pa_{\phi_i}f^\mathrm{loc} - \beta_i\Delta\phi_i, \quad i \in \mathcal{S}^*, & \mbox{in }Q_T, \\
		\label{eq.Phiv}
		&\pa_t\Phi_v + \vel\cdot\nabla\Phi_v = 
		-\gamma_v(\phi_{\mathcal{S}^*})\left(
		\pa_{\Phi_v}f^\mathrm{loc} -\beta_v \Delta\Phi_v
		\right) & \mbox{in }Q_T.
	\end{align}
    The mobilities $\gamma_{ij} = \gamma_{ij}(\phi_{\mathcal{S}^*},\Phi_v)$ form a symmetric, positive semidefinite matrix. Specifically, we set
	\begin{align*}
		\gamma_{ij} &= (1-\omega(\Phi_v))\gamma_{ij}^\mathrm{liq} + \omega(\Phi_v)\gamma_{ij}^\mathrm{vap}, 
	\end{align*}
	where $\gamma_{ij}^{\mathrm{liq}}(\phi_{\mathcal{S}^*})$ and $\gamma_{ij}^\mathrm{vap}(\phi_{\mathcal{S}^*})$ are symmetric, positive semidefinite matrices, depending only on the volume fractions $\phi_{\mathcal{S}^*}$. The interpolating function $\omega\in C^2(\R)$ satisfies $\omega=0$ on $(-\infty,0]$, $\omega=1$ on $[1,\infty)$, and is increasing on $(0,1)$. The columns of the mobility matrices (and hence, by symmetry, the rows) sum to zero, $$ \sum_{j \in \mathcal{S}^*}\gamma_{ij}^\mathrm{liq} = \sum_{j \in \mathcal{S}^*} \gamma_{ij}^\mathrm{vap} = 0,\quad i \in \mathcal{S}^*. $$
    The fluxes \eqref{eq.Ji} can be rewritten as
	\begin{align*}
		\bJ_i = 
        -\sum_{j \in \mathcal{S}} \gamma_{ij}\nabla(\mu_j-\mu_a), \ i \in \mathcal{S}, \ \text{ and } \ 
		\bJ_a = -\sum_{i \in \mathcal{S}} \bJ_i ,
	\end{align*}
	and, since $\phi_a = 1-\sum_{i \in \mathcal{S}} \phi_i$, the equation for $\phi_a$ in \eqref{eq.phii} can be dropped. 
    The (local part of the) nondimensional free energy density is
	\begin{align}\label{f}
		f^{\mathrm{loc}} = &\sum_{i \in \mathcal{S}^*} \frac{\phi_j}{N_j}\log\phi_j + (1-\omega(\Phi_v))\sum_{i,j\in \mathcal{S}^*} \chi_{ij}\phi_i\phi_j -
		\omega(\Phi_v)\sum_{i \in \mathcal{S}^*} \frac{\phi_j}{N_j}\log\phi_j^\mathrm{sat}. 
	\end{align}
    The system is subject to the boundary conditions
    \begin{align*}
        \sum_{j \in \mathcal{S}^*} \gamma_{ij} \nabla \mu_j \cdot \nu  &= 0 \quad \text{on} \quad \partial \Omega \times (0,T), \ i \in \mathcal{S}^*, \\
        \partial_\nu \phi_i &= 0 \quad \text{on} \quad \partial \Omega \times (0,T),\ i \in \mathcal{S}^*, \\
        \partial_\nu \Phi_v &= 0 \quad \text{on} \quad \partial \Omega \times (0,T), \\ 
        \vel \cdot \nu &= 0 \quad \text{on} \quad \partial \Omega \times (0,T), \\ 
        (\bI-\nu\otimes\nu)D(\vel)\nu &= 0  \quad \text{on} \quad \partial \Omega \times (0,T). 
    \end{align*}
    \section{Analysis} \label{sec.analysis}
    In this section, we prove the global-in-time existence of weak solutions to the PDE system presented in Subsection~\ref{sec.nondimeq}, in a slightly more general setting. For readability, we write the volume fractions as $\phi_1, \dots, \phi_n$. The model of Subsection~\ref{sec.nondimeq} corresponds to the case $n=4$, identifying $\mathcal{S}^*$ with $\{ 1, \dots, n\}$. Furthermore, we redefine the notation
	\begin{align*}
		\vec\phi = (\phi_1,\ldots,\phi_{n-1}) \  \text{ and } \ \phi_n = 1-\sum_{i=1}^{n-1}\phi_i,
	\end{align*}
    and define 
    \begin{align*}
        \mathscr{D}=\left\{
		X\in [0,1]^{n-1}~:~~\sum_{i=1}^{n-1}X_i\leq 1
		\right\} \  \text{ and } \  \mathscr{S} = \left\{
        Y\in [0,1]^n ~:~ \sum_{i=1}^{n}Y_i=1\right\}.
    \end{align*}
    The assumptions of nonnegativity of all phases and the volume constraint can be restated as $\vec\phi(Q_T)\subset \mathscr{D}$.
    Moreover, we denote by $\mathfrak{q}_i \equiv\mu_i-\mu_n$ the relative chemical potentials. Lastly, we couple \eqref{eq.v}--\eqref{eq.Phiv} with the following boundary conditions
    \begin{align}\label{bc}
		\begin{cases}
			\vel\cdot\nu = 0,\quad (I-\nu\otimes\nu)D(\vel)\nu = 0,\\
			\pa_\nu\phi_i = 0,\quad  
			\bJ_i\cdot\nu = 0,\quad 
			\pa_\nu\Phi_v = 0,
		\end{cases}
		\qquad\mbox{on }\pa\Omega\times (0,T)\quad 
		(i=1,\ldots,n)
	\end{align}
    and generic initial conditions
	\begin{align}\label{ic}
		\vel(\cdot,0)=\vel^{\mathrm{in}},\quad\phi_i(\cdot,0)=\phi_{i}^{\mathrm{in}},\quad
		\Phi_v(\cdot,0) = \Phi_{v}^{\mathrm{in}}
		\qquad\mbox{in }\Omega\quad (i=1,\ldots,n).
	\end{align}
    The initial data satisfy
	\begin{align}
    \label{hp.init}
    &\vel^{\mathrm{in}}\in L_\sigma^2(\Omega),\quad 
    \phi_i^{\mathrm{in}}\in H^1(\Omega),\quad
    \phi_i^{\mathrm{in}}>0,\quad\sum_{j=1}^n\phi_j^{\mathrm{in}}=1\quad\mbox{a.e.~in }\Omega,\\
    \label{hp.init.2}
    &\Phi_v^{\mathrm{in}}\in H^1(\Omega),\quad 0\leq \Phi_v^{\mathrm{in}}\leq 1\quad\mbox{a.e.~in }\Omega,
    \end{align}
    for $i=1,\ldots,n$, where $L^2_\sigma(\Omega)$ is defined in Subsection~\ref{sec.weakform}.
    
    \subsection{Preliminary considerations}\label{ss.ana.0}
    In this subsection, we state assumptions on the functions appearing in the system presented in Subsection~\ref{sec.nondimeq}, introduce an energy functional, and show that it satisfies an evolution equation.
    
    \paragraph{\bf Assumptions.}
    We assume that the mobilities for the liquid and vapor phases satisfy the following assumptions
    for a suitable choice of 
	$D_1,\ldots,D_n,c_\gamma>0$, $0<\epsilon_0<\frac{1}{2}$:
	\begin{align}
		\label{hp.mob.1}
		\left|\gamma_{ij}^{L,V}\right|
		+\phi_\ell^{\epsilon_0}\,\left|\pa_{\phi_\ell}\gamma_{ij}^{L,V}\right|\leq C\phi_j^{\epsilon_0}\quad 
		\forall (\phi_1,\ldots,\phi_n)\in\mathscr{S},\\
		\label{hp.mob.2}
		\sum_{k,s=1}^n\gamma_{ks}^{L,V}\xi_k\xi_s\geq c_\gamma
		\sum_{k,s=1}^n\left(
		D_k\phi_k\delta_{ks} - \frac{D_k D_s\phi_k\phi_s}{\sum_{m=1}^n D_m\phi_m}
		\right)\xi_k\xi_s,\\ \nonumber
		\forall\xi\in\R^n,~~ 
		\forall (\phi_1,\ldots,\phi_n)\in\mathscr{S},\\
		\label{hp.mob.3}
		(\phi_1,\ldots,\phi_n)\in\mathscr{S}\mapsto\gamma_{ij}^{L,V}(\phi_1,\ldots,\phi_n)\in\R\quad\mbox{is of class }C^1(\mathscr{S}),
	\end{align}
	for $i,j,\ell=1,\ldots,n$.\\
    The mobility function $\gamma_v:\dom\to\R$ is assumed to be $C^{1}(\dom)$, while the viscosity function $\eta:\dom\to\R$ is assumed to be Lebesgue measurable. Furthermore, we assume both functions to be uniformly positive definite and bounded:
	\begin{align}
		\label{hp.gammaV}
		0<\gamma_{v,0}\leq \gamma_v(\vec\phi)\leq \gamma_{v,1}\quad\forall\vec\phi\in\dom ,\\
		\label{hp.eta}0<\eta_0\leq\eta(\vec\phi)\leq\eta_1\quad\forall\vec\phi\in\dom .
	\end{align}
	
	\begin{rem}
		Assumptions \eqref{hp.mob.1}--\eqref{hp.mob.3} are mainly needed to obtain a priori estimates for the liquid phases, the chemical potentials, and the diffusion fluxes, as well as ensuring a well-posed weak formulation. The lower bound \eqref{hp.mob.2} is inspired by the expression for the mobilities in the so-called ``slow-mode theory''~\cite{RonHar22}.
	\end{rem}

	\paragraph{\bf Energy balance.}
    Motivated by \eqref{modeling.dissip.ineq}, we define the {\em (mathematical) energy density}
	$$
	H := f + \frac{1}{2}\rho |\vel|^2.
	$$
    More explicitly,
	\begin{align}\label{E}
		H =& \frac{1}{2}\rho |\vel|^2
		+ \frac{1}{2}\sum_{j=1}^n\beta_j|\nabla\phi_j|^2
		+ \frac{1}{2}\beta_v|\nabla\Phi_v|^2\\ \nonumber
		&+ \sum_{j=1}^n\frac{\phi_j}{N_j}\log\phi_j + (1-\omega(\Phi_v))\sum_{i,j=1}^n\chi_{ij}\phi_i\phi_j - 
		\omega(\Phi_v)\sum_{j=1}^n\frac{\phi_j}{N_j}\log{\phi_j^{sat}}.
	\end{align}
    The following {\em energy balance} is derived in Subsection~\ref{sec.thdyn} and Appendix~\ref{app.thdyncons} via physical and thermodynamic arguments. Here, we present an alternative derivation, which will be useful in later sections.  
	\begin{lemma}\label{lem.energy}
		The {\em energy balance} holds:
		\begin{align}\label{en.bal}
			\pa_t H + \Div \mathcal{J}_H + \Xi_H = 0, 
		\end{align}
		where the {\em energy flux} $\mathcal{J}_H$ and the {\em energy dissipation} $\Xi_H$ are given by
		\begin{align*}
			\mathcal{J}_H =& H\vel 
            +\sum_{i=1}^{n}\left(
			\frac{1}{2}|\vel|^2\rhob_i + \mu_i\right)\bJ_i 
			- 2\eta D(\vel)\vel + p\vel
            -\sum_{j=1}^{n}\beta_j\pa_t\phi_j\nabla\phi_j 
			-\beta_v\pa_t\Phi_v\nabla\Phi_v ,\\
			\Xi_H =& 2\eta|D(\vel)|^2 + \sum_{i,j=1}^{n-1}\gamma_{ij}\nabla(\mu_i - \mu_n)\cdot\nabla(\mu_j - \mu_n) 
			+ \gamma_v|\pa_{\Phi_v}f - \beta_v\Delta\Phi_v|^2. 
		\end{align*}
	\end{lemma}
	\begin{proof}
		We multiply \eqref{eq.v} by $\vel$ and employ the incompressibility constraint \eqref{eq.divv} to get
		\begin{align}\label{rDv2}
			\rho &\pa_t\left(\frac{1}{2}|\vel|^2\right) + 
			\left(
			\rho\vel - \sum_{i,j=1}^{n-1}(\rhob_i - \rhob_n)\gamma_{ij}\nabla(\mu_j - \mu_n)
			\right)\cdot\nabla\frac{|\vel|^2}{2}
			+2\eta|D(\vel)|^2\\ \nonumber
			&= 
			\Div(2\eta D(\vel)\vel - p\vel) - \sum_{j=1}^{n}\beta_j\vel\cdot\Div(\nabla\phi_j\otimes\nabla\phi_j) - \beta_v\vel\cdot\Div(\nabla\Phi_v\otimes\nabla\Phi_v).
		\end{align}
		We multiply \eqref{eq.phii} by $\rhob_i$, sum over $i=1,\ldots,n$, and using $\sum_{i=1}^n \bJ_i = 0 $ to deduce the total mass balance
		\begin{align*}
			\pa_t\rho + \vel\cdot\nabla\rho = \Div\left(
			\sum_{i,j =1}^{n-1}(\rhob_i - \rhob_n)\gamma_{ij}\nabla(\mu_j - \mu_n)
			\right).
		\end{align*}
		Hence, multiplying the above equation by $|\vel|^2/2$ and adding the resulting identity to \eqref{rDv2}, we obtain the kinetic energy balance
		\begin{align}\label{kin.bal}
			\pa_t&\left(\frac{1}{2}\rho|\vel|^2\right) + 
			\vel\cdot\nabla\left(\frac{1}{2}\rho|\vel|^2\right)
			+\Div\left(
            p\vel
			- \frac{|\vel|^2}{2}\sum_{i,j=1}^{n-1}(\rhob_i - \rhob_n)\gamma_{ij}\nabla(\mu_j - \mu_n)
			- 2\eta D(\vel)\vel 
			\right) \\ \nonumber 
			&= 
			-\sum_{j=1}^{n}\beta_j\vel\cdot\Div(\nabla\phi_j\otimes\nabla\phi_j) - \beta_v\vel\cdot\Div(\nabla\Phi_v\otimes\nabla\Phi_v)
			-2\eta|D(\vel)|^2.
		\end{align}
        Recall that $\mu_i = \pa_{\phi_i}f - \beta_i\Delta\phi_i$.
		We multiply \eqref{eq.phii} by $\mu_i$ and sum over $i=1,\ldots,n$, obtaining
		\begin{align*}
			\sum_{i=1}^{n} &\pa_{\phi_i}f\,\pa_t\phi_i 
			+ \pa_t\sum_{i=1}^{n}\frac{\beta_i}{2}|\nabla\phi_i|^2
			+ \vel\cdot\left[
			\sum_{i=1}^{n}\pa_{\phi_i}f\,\nabla\phi_i 
			+\nabla \sum_{i=1}^{n}\frac{\beta_i}{2}|\nabla\phi_i|^2
			\right]\\ \nonumber
			=&\Div\left(
			\sum_{i=1}^{n}\beta_i\pa_t\phi_i \nabla\phi_i
			+\sum_{i,j=1}^{n-1}(\mu_i - \mu_n)\gamma_{ij}\nabla(\mu_j - \mu_n)
			\right) + \sum_{j=1}^{n}\beta_j\vel\cdot\Div(\nabla\phi_j\otimes\nabla\phi_j)\\ \nonumber
			& -\sum_{i,j=1}^{n-1}\gamma_{ij}\nabla(\mu_i - \mu_n)\cdot\nabla(\mu_j - \mu_n).
		\end{align*}
		By multiplying \eqref{eq.Phiv} by $\pa_{\Phi_v}f - \beta_v\Delta\Phi_v$ and performing similar computations, we also obtain
		\begin{align}\label{df.Phiv}
			\pa_{\Phi_v}f\,&\pa_t\Phi_v 
			+ \pa_t \frac{\beta_v}{2}|\nabla\Phi_v|^2
			+ \vel\cdot\left[
			\pa_{\Phi_v}f\,\nabla\Phi_v 
			+\nabla\frac{\beta_v}{2}|\nabla\Phi_v|^2
			\right]\\ \nonumber
			=&\Div\left(
			\beta_v\pa_t\Phi_v \nabla\Phi_v\right) + \beta_v\vel\cdot\Div(\nabla\Phi_v\otimes\nabla\Phi_v) -\gamma_v|\pa_{\Phi_v}f - \beta_v\Delta\Phi_v|^2.
		\end{align}
		Summing \eqref{kin.bal}--\eqref{df.Phiv} yields the statement. This finishes the proof of the lemma.
	\end{proof}

	\subsection{A priori estimates}\label{sec.apriori-est}
	Here, we derive some a priori estimates for the solutions to \eqref{eq.v}--\eqref{eq.Phiv}, which will be employed to derive the weak formulation of the system.
	
	Integrating the energy equation \eqref{en.bal} over space and time yields the following bounds:
	\begin{align}\label{est.LinfL2}
	\sum_{j=1}^{n-1}\|\nabla\phi_j\|_{L^\infty(0,T; L^2(\Omega))} + 
		\|\nabla\Phi_v\|_{L^\infty(0,T; L^2(\Omega))} +
		\|\vel\|_{L^\infty(0,T; L^2(\Omega))} &\leq C,\\
		\label{est.L2L2}
		\|D(\vel)\|_{L^2(0,T; L^2(\Omega))}  
		+ \|\gamma_v^{1/2}
		\Delta\Phi_v\|_{L^2(0,T; L^2(\Omega))} &\leq C,\\
		\label{est.namu}
		\int_0^T\int_\Omega \sum_{i,j=1}^{n-1}\gamma_{ij}(\vec\phi,\Phi_v)\nabla(\mu_i-\mu_n)\cdot\nabla(\mu_j-\mu_n)\dx\dt &\leq C.
	\end{align}	
    The following lemma contains crucial estimates for the derivatives of the phase fields $\phi_i$, $i=1,\ldots,n$ and is derived from \eqref{est.namu}.

	\begin{lemma}\label{lem.est.namu}
		The following a priori estimates hold for every $i=1,\ldots,n$:
		\begin{align}\label{bound.sqrtphii}
			\|\nabla\sqrt{\phi_i}\|_{L^2(0,T; L^2(\Omega))} &\leq C,\\
			\label{bound.Deltaphii}
			\|D^2{\phi_i}\|_{L^2(0,T; L^2(\Omega))} &\leq C,\\
			\label{bound.L.mu}
			\int_0^T\int_\Omega\sum_{k=1}^{n}\phi_k\left|
			\nabla\mu_k - \sum_{j=1}^{n}\alpha_j\nabla\mu_j
			\right|^2\dx\dt &\leq C,
			\qquad\alpha_i \equiv \frac{D_i\phi_i}{\sum_{k=1}^n D_k\phi_k}.%
		\end{align}
		In particular, $\phi_i^\epsilon\mu_i$ is bounded in $L^2(0,T; L^2(\Omega))$ for every $\epsilon>0$. Furthermore, it holds that
		\begin{align}\label{bound.flux}
			\|J^{V,i}\|_{L^2(0,T; L^2(\Omega))} + 
			\|J^{L,i}\|_{L^2(0,T; L^2(\Omega))}\ &\leq C,
		\end{align}
		where $J^{V,i}$, $J^{L,i}$ are the vapor and liquid components of the diffusion flux, i.e.
		\begin{align*}
			J^{V,i} &= -\omega(\Phi_v)\sum_{j=1}^{n-1}\gamma_{ij}^V\nabla(\mu_j-\mu_n),\\
			J^{L,i} &= -(1-\omega(\Phi_v))\sum_{j=1}^{n-1}\gamma_{ij}^L\nabla(\mu_j-\mu_n).
		\end{align*}
		
	\end{lemma}
	\begin{proof}
		Ineq.~\eqref{bound.L.mu} is a straightforward consequence of eq.~\eqref{est.namu} and eq.~\eqref{hp.mob.2}. Eq.~\eqref{eq.mui} and eq.~\eqref{bound.L.mu} imply
		\begin{align}\label{bound.prefactor}
			\int_0^T\int_\Omega\sum_{i=1}^{n}
			N_i\phi_i\left|
			\nabla\pa_{\phi_i}f - \sum_{j=1}^{n}\alpha_j\nabla\pa_{\phi_j}f
			-\left(
			\beta_i\nabla\Delta\phi_i -
			\sum_{j=1}^n\alpha_j\beta_j\nabla\Delta\phi_j
			\right)
			\right|^2\dx\dt\leq C .
		\end{align}
		Computing the square inside the integral yields
		\begin{align}\label{est.namu.5} 
			A &+ B + c\Gamma \leq C, \text{ where }\\
			A &= \int_{0}^{T}\int_\Omega\sum_{i=1}^n\phi_i \left|
			\nabla\pa_{\phi_i}f - \sum_{j=1}^{n}\alpha_j\nabla\pa_{\phi_j}f
			\right|^2\dx\dt , \notag\\ \label{est.namu.B}
			B &= \int_{0}^{T}\int_\Omega\sum_{i=1}^n\phi_i \left|
			\beta_i\nabla\Delta\phi_i -
			\sum_{j=1}^n\alpha_j\beta_j\nabla\Delta\phi_j
			\right|^2\dx\dt ,\\ \label{est.namu.Gamma}
			\Gamma &= -\int_{0}^{T}\int_\Omega\sum_{i=1}^n
			N_i\phi_i \left(
			\nabla\pa_{\phi_i}f - \sum_{j=1}^{n}\alpha_j\nabla\pa_{\phi_j}f
			\right)\\ \nonumber
			&\qquad \cdot\left(
			\beta_i\nabla\Delta\phi_i -
			\sum_{j=1}^n\alpha_j\beta_j\nabla\Delta\phi_j 
			\right)\dx\dt .
		\end{align}
		We now compute
		\begin{align*}
			\sum_{j=1}^{n}\alpha_j\nabla\pa_{\phi_j}f = &
			\sum_{j=1}^{n}\alpha_j\nabla\left(
			\frac{\log\phi_j + 1}{N_j} + 2(1-\omega(\Phi_v))\sum_{k=1}^{n}\chi_{jk}\phi_k - \omega(\Phi_v)\frac{1}{N_j}\log{\phi_j^{sat}} 
			\right)\\
			=&\frac{\sum_{j=1}^n D_j N_j^{-1}\nabla\phi_j}{\sum_{k=1}^{n}D_k\phi_k} +
			2(1-\omega(\Phi_v))\frac{\sum_{j,k=1}^{n}\chi_{jk}D_j\phi_j\nabla\phi_k}{\sum_{k=1}^{n}D_k\phi_k} \\ 
            &- \frac{\sum_{j=1}^n D_j \phi_j N_j^{-1}\log\phi_j^{sat} + 2\sum_{j,k=1}^{n}\chi_{jk}D_j\phi_j\phi_k}{\sum_{k=1}^{n}D_k\phi_k}\nabla\omega(\Phi_v).
		\end{align*}
		According to \eqref{est.LinfL2}, 
		\begin{align*}
			\left\| \sum_{j=1}^{n}\alpha_j\nabla\pa_{\phi_j}f\right\|_{
				L^\infty(0,T; L^2(\Omega))}\leq C .
		\end{align*}
		Furthermore, 
		\begin{align*}
			\nabla\pa_{\phi_i}f = N_i^{-1}\phi_i^{-1}\nabla\phi_i + \nabla\left[
			2(1-\omega(\Phi_v))\sum_{k=1}^{n}\chi_{ik}\phi_k
			\right] - \left(\frac{1}{N_i}\log\phi_i^{sat} \right)\nabla \omega(\Phi_v),
		\end{align*}
		and hence 
		\begin{align}\label{namu.decomp}
			&\nabla\pa_{\phi_i}f - \sum_{j=1}^{n}\alpha_j\nabla\pa_{\phi_j}f
			= N_i^{-1}\phi_i^{-1}\nabla\phi_i + g_i,\\ 
			&g_i \equiv \nabla\left[
			2(1-\omega(\Phi_v))\sum_{k=1}^{n}\chi_{ik}\phi_k
			\right] - \left(\frac{1}{N_i}\log\phi_i^{sat} \right)\nabla \omega(\Phi_v)
			- \sum_{j=1}^{n}\alpha_j\nabla\pa_{\phi_j}f, \notag\\
			\label{namu.decomp.3}
			&\|g_i\|_{L^\infty(0,T; L^2(\Omega))}\leq C .
		\end{align}
		It follows that 
		\begin{align*}
			A &\geq \frac{1}{2}\int_0^T\int_\Omega\sum_{i=1}^{n}
            N_i^{-2}\phi_i|\phi_i^{-1}\nabla\phi_i|^2\dx\dt 
			-\int_0^T\int_\Omega\sum_{i=1}^{n}\phi_i|g_i|^2\dx\dt ,
		\end{align*}
		resulting in 
		\begin{align}\label{A.lb}
			A &\gtrsim \int_0^T\int_\Omega\sum_{i=1}^{n}|\nabla\sqrt{\phi_i}|^2 
			\dx\dt - C.
		\end{align}
		Next, we find an expression for $\Gamma$. From \eqref{est.namu.Gamma} and \eqref{namu.decomp} we obtain
		\begin{align*}
			&\Gamma = -\int_0^T\int_\Omega\sum_{i=1}^{n}
			\left(
			\nabla\phi_i + N_i \phi_i g_i
			\right)\cdot\left(
			\beta_i\nabla\Delta\phi_i -
			\sum_{j=1}^n\alpha_j\beta_j\nabla\Delta\phi_j
			\right)\dx\dt = \Gamma_1 + \Gamma_2,\\
			&\Gamma_1 = -\int_0^T\int_\Omega\sum_{i=1}^{n}\nabla\phi_i
			\cdot\left(
			\beta_i\nabla\Delta\phi_i -
			\sum_{j=1}^n\alpha_j\beta_j\nabla\Delta\phi_j
			\right)\dx\dt,\\
			&\Gamma_2 = -\int_0^T\int_\Omega\sum_{i=1}^{n}
			N_i \phi_i g_i
			\cdot\left(
			\beta_i\nabla\Delta\phi_i -
			\sum_{j=1}^n\alpha_j\beta_j\nabla\Delta\phi_j
			\right)\dx\dt .
		\end{align*}
		Given that the $\phi_i$ sum to 1, we immediately obtain
		\begin{align}\label{Gamma1}
			\Gamma_1 = -\int_{0}^T\int_\Omega\sum_{i=1}^n\beta_i\nabla\phi_i\cdot\nabla\Delta\phi_i\dx\dt = \int_0^T\int_\Omega\sum_{i=1}^n\beta_i(\Delta\phi_i)^2\dx\dt ,
		\end{align}
		while for arbitrary $\epsilon>0$, using Young's inequality, eq.~\eqref{est.namu.B}, and eq.~\eqref{namu.decomp.3},
		\begin{align*}
			\Gamma_2\geq -\epsilon\int_0^T\int_\Omega\sum_{i=1}^{n}
			\phi_i\left|
			\beta_i\nabla\Delta\phi_i -
			\sum_{j=1}^n\alpha_j\beta_j\nabla\Delta\phi_j
			\right|^2\dx\dt - C(\epsilon) = -\epsilon B - C(\epsilon).
		\end{align*}
		Consequently,  
		\begin{align}\label{Gamma.lb}
			\Gamma\geq \int_0^T\int_\Omega\sum_{i=1}^n\beta_i(\Delta\phi_i)^2\dx\dt 
			-\epsilon B - C(\epsilon) .
		\end{align}
		Combining eq.~\eqref{est.namu.5}, eq.~\eqref{A.lb}, eq.~\eqref{Gamma.lb}, and choosing $\epsilon>0$ small enough yields
		\begin{align*}
			\int_0^T\int_\Omega &\sum_{i=1}^{n}\left[
			|\nabla\sqrt{\phi_i}|^2 + (\Delta\phi_i)^2
			\right]\dx\dt + \int_{0}^{T}\int_\Omega\sum_{i=1}^n\phi_i \left|
			\beta_i\nabla\Delta\phi_i -
			\sum_{j=1}^n\alpha_j\beta_j\nabla\Delta\phi_j
			\right|^2\dx\dt\leq C .
		\end{align*}
		Hence, eq.~\eqref{bound.sqrtphii} is proven, as well as 
		\begin{align*}
			\|\Delta\phi_i\|_{L^2(Q_T)}\leq C, \quad i=1,\ldots,n,
		\end{align*}
		which, together with the boundary condition~\eqref{bc} for $\phi_i$, yields \eqref{bound.Deltaphii} via standard elliptic regularity arguments.

		Bound \eqref{bound.flux} follows from \eqref{est.namu}, given that the matrices
		$\gamma_{ij}^L$, $\gamma_{ij}^V$ are symmetric, positive semidefinite, and bounded.
	\end{proof}
    \begin{rem}
    The proof of \eqref{bound.sqrtphii} and \eqref{bound.Deltaphii} bears some similarities with the proof of the same bounds in~\cite[Lemma 4.4]{EhrMarPie21}. While the two strategies differ due to the different structure of the mobility matrices, in both cases a cancellation property, derived from the fact that the phases sum to 1, is crucial to the argument. In Lemma~\ref{lem.est.namu}, this property is apparent in \eqref{Gamma1}. We also point out that the form of the term $\Gamma_1$ as it appears in \eqref{Gamma1} stems from the specific choice of the prefactor $N_i^{-1}$ in \eqref{bound.prefactor}.
    \end{rem}
    
	The following Gagliardo--Nirenberg--Sobolev inequality holds:
	\begin{align}
		\label{GNS.10.3}
		\|u\|_{L^{10/3}(Q_T)}\leq C \|u\|_{L^2(0,T; H^1(\Omega))}^{3/5}\|u\|_{L^\infty(0,T; L^2(\Omega))}^{2/5}
	\end{align}
    for every $u\in L^2(0,T; H^1(\Omega))\cap L^\infty(0,T; L^2(\Omega))$. Hence
    \eqref{est.LinfL2}, \eqref{est.L2L2}, \eqref{bound.Deltaphii}, \eqref{GNS.10.3} imply
	\begin{align}
		\label{est.Lp}
		\|\nabla\phi_i\|_{L^{10/3}(Q_T)}
		+\|\nabla\Phi_v\|_{L^{10/3}(Q_T)}
		+\|\vel\|_{L^{10/3}(Q_T)}
		\leq C,\quad i=1,...,n.
	\end{align}

	\subsection{Weak formulation of the system}
	\label{sec.weakform}
	We define the function spaces
	\begin{align*}
		L^p_\sigma(\Omega,\R^3)\equiv &\left\{
		\bm{w}\in L^p(\Omega,\R^3)~:~\int_\Omega\bm{w}\cdot\nabla\varphi
		\dx = 0\quad\forall\varphi\in W^{1,p'}(\Omega)
		\right\},\\
		W^{1,p}_\sigma(\Omega,\R^3)\equiv &\left\{
		\bm{w}\in W^{1,p}(\Omega,\R^3)~:~\Div\bm{w}=0~\mbox{a.e.~in }\Omega,~\bm{w}\cdot\nu =0~\mbox{on }\pa\Omega\right\},\\
		H^1_\sigma(\Omega,\R^3)\equiv & W^{1,2}_\sigma(\Omega,\R^3),
	\end{align*}
    and the function
    \begin{align}\label{Q}
    \mathcal{Q}(\vec\phi,\Phi_v)[\bJ_1,\ldots,\bJ_n] = \sup\left\{
    2\sum_{i=1}^n \bJ_i\cdot\bm{\zeta}_i - \sum_{i,j=1}^n\gamma_{ij}(\vec\phi,\Phi_v)\bm{\zeta}_i\cdot\bm{\zeta}_j~\Bigg\vert~\bm{\zeta}_1,\ldots,\bm{\zeta}_n\in\R^3
    \right\},
    \end{align}
    which is simply the Legendre transform of the quadratic form associated with $\gamma_{ij}$. Since $(\gamma_{ij})_{i,j=1,\ldots,n}$ is symmetric and positive definite on $\mbox{Span}\{ {\bf 1} \}^\perp$ for fixed $\vec\phi\in\mathring{\dom}$, $\Phi_v\in \R$, it is easy to see that $\mathcal{Q}(\vec\phi,\Phi_v)[\bJ_1,\ldots,\bJ_n]$ coincides with $\sum_{i,j=1}^n\gamma_{ij}\nabla\mu_i\cdot\nabla\mu_j$ when $\vec\phi\in\mathring{\dom}$, while for $\vec\phi\in\dom$ it is a convex and well-defined mapping $\R^{3n}\to\R\cup\{+\infty\}$. For these reasons, we employ $\mathcal{Q}(\vec\phi,\Phi_v)[\bJ_1,\ldots,\bJ_n]$ in the energy balance instead of $\sum_{i,j=1}^n\gamma_{ij}\nabla\mu_i\cdot\nabla\mu_j$.
	
	\begin{defn}[Weak solution to \eqref{eq.v}--\eqref{eq.Phiv}]\label{def.weaksol} 
		A measurable function $(\vel,\mu,\phi,\Phi_v) : Q_T\to \R^3\times\R^n\times\R^n\times\R$ is called a {\em weak solution} to \eqref{eq.v}--\eqref{eq.Phiv}, \eqref{bc}, \eqref{ic} if
		\begin{align*}
			\phi_i \geq 0\quad i=1,\ldots,n,\qquad\sum_{i=1}^n\phi_i=1,\qquad
			0\leq\Phi_v\leq 1\qquad\mbox{a.e.~in }Q_T,
		\end{align*}
		the functions have the regularity
		\begin{align*}
			&\vel\in L^2(0,T; H_\sigma^1(\Omega))\cap L^\infty(0,T; L^2(\Omega)),\quad  
            \pa_t(\rho\vel)\in L^{5/4}(0,T; (W_\sigma^{1,5}(\Omega))'),\\
			&\phi_i\in L^2(0,T; H^2(\Omega))
			\cap L^\infty(0,T; H^1(\Omega))
			\cap H^1(0,T; (H^1(\Omega))'),\\
			&\sqrt{\phi_i}\in L^2(0,T; H^1(\Omega)),\quad
			\phi_i^\epsilon\mu_i\in L^2(Q_T),\\
			&\Phi_v\in L^2(0,T; H^2(\Omega))
			\cap L^\infty(0,T; H^1(\Omega))
			\cap W^{1,\frac{5}{3}}(0,T; L^{\frac{5}{3}}(\Omega)),
		\end{align*}
		for $i=1,\ldots,n$ and for every $\epsilon>0$,
		and the following relations hold:
		\begin{enumerate}
\item {\em Weak formulation of the momentum equation:}
\begin{align}\label{weak.mom}
\int_0^T\langle \pa_t(\rho\vel), \bm{w}\rangle_{((W_\sigma^{1,5}(\Omega))',W_\sigma^{1,5}(\Omega))}\dt 
=\int_0^T\int_\Omega\left(
\rho\vel\otimes\vel 
+ \vel\otimes\sum_{i=1}^{n}\rhob_i \bJ_i
-\mathbb{S}\right):\nabla\bm{w}\dx\dt\\ \nonumber
\forall \bm{w}\in L^5(0,T; W_\sigma^{1,5}(\Omega)),
\end{align}
with the viscous stress $\mathbb{S}$ given by
\begin{align}\label{weak.mom.M}
\mathbb{S} = 2\eta(\vec\phi)D(\vel)
-\sum_{j=1}^{n}\beta_j \nabla\phi_j\otimes\nabla\phi_j
-\beta_v \nabla\Phi_v\otimes\nabla\Phi_v .
\end{align}

\item {\em Weak formulation of the liquid-phase equations:}
			\begin{align}
				\label{weak.phii}
				\int_0^T &\sum_{i=1}^{n}
				\langle\pa_t\phi_i,\psi_i\rangle_{((H^1(\Omega))',H^1(\Omega))}\dt
				=\int_0^T\int_\Omega\sum_{i=1}^{n}
				\nabla\psi_i\cdot (\bJ_i + \phi_i\vel) \dx\dt 
				\\ \nonumber
				&\forall\psi_1,\ldots,\psi_{n}
				\in L^2(0,T; H^1(\Omega)),
			\end{align}
			where the diffusive fluxes $\bJ_1,\ldots,\bJ_n\in L^2(Q_T; \R^3)$ satisfy
			\begin{align}
				\label{Ji}
				\sum_{i=1}^n &\int_0^T\int_\Omega\bm{\Psi}_i\cdot \bJ_i\dx\dt 
				= 
				\sum_{i,j=1}^n\int_0^T \int_\Omega\gamma_{ij}\mu_j\,\Div\bm{\Psi}_i\dx\dt\\ \nonumber
				&+\sum_{i,j=1}^n\int_0^T \int_\Omega\bm{\Psi}_i\cdot
				\left(
				\sum_{\ell=1}^n\frac{\phi_\ell^{\epsilon_0}}{1-\epsilon_0}\,\pa_{\phi_\ell}\gamma_{ij}\nabla\phi_\ell^{1-\epsilon_0}
				+\pa_{\Phi_v}\gamma_{ij}\nabla\Phi_v
				\right)\mu_j\dx\dt\\ \nonumber
				\forall & \bm{\Psi}_1,\ldots,\bm{\Psi}_n\in L^2(0,T; H_0^1(\Omega;\R^3))\cap L^\infty(Q_T;\R^3).
			\end{align}
			\item {\em Definition of the chemical potentials:}
			Eq.~\eqref{eq.mui} holds a.e.~in $Q_T\cap\{\phi_i>0\}$, i.e. 
			\begin{align}\label{weak.mui}
				\begin{cases}
					\mu_i = \pa_{\phi_i}f^{loc}(\phi_1,\ldots,\phi_n,\Phi_v) - \beta_i\Delta\phi_i \quad \mbox{a.e.~in }Q_T\cap\{\phi_i>0\}, \\ \smallskip 
					\pa_\nu\phi_i = 0\quad\mbox{in }L^2(0,T; H^{1/2}(\pa\Omega)),
				\end{cases}
			\end{align}
			for $i=1,\ldots,n$.
			\item {\em Equation for the vapor phase:}
			Eq.~\eqref{eq.Phiv} holds a.e.~in $Q_T$ and \\ $\pa_\nu\Phi_v = 0$ in $L^2(0,T; H^{1/2}(\pa\Omega))$.
	\item {\em Energy balance:} for almost every $0\leq t_1<t_2\leq T$ it holds
	\begin{align}\label{weak.en.bal}
		&\int_\Omega H(t_2)\dx + \int_{t_1}^{t_2}
		\int_\Omega\Xi_H\dx\dt
		\leq \int_\Omega H(t_1)\dx 
	\end{align}
	with the energy $H$ as in \eqref{E} and the energy dissipation $\Xi_H$ given by
    \begin{align}
		\label{weak.en.bal.Xi}
		\Xi_H =&
		2\eta(\vec\phi)|D\vel|^2 
        +\mathcal{Q}(\bJ_1,\ldots,\bJ_n)
		+\gamma_v(\vec\phi)
		\left|\pa_{\Phi_v}f^{loc}(\vec\phi,\Phi_v)
		-\beta_v\Delta\Phi_v\right|^2 .
	\end{align} 
    \item {\em Initial conditions:} the following convergence relations hold for $t\to 0^+$, $i=1,\ldots,n$:
    \begin{align}
    \label{ic.weak}
    &\vel(t)\to \vel^{\mathrm{in}}\quad\mbox{strongly in }L^{2}(\Omega),\quad 
    \phi_i(t)\to\phi_i^{\mathrm{in}},\quad\Phi_v(t)\to \Phi_v^{\mathrm{in}}\quad\mbox{strongly in }H^1(\Omega).
    \end{align}
    \end{enumerate}
	\end{defn}

    \begin{rem}
    We point out that \eqref{weak.mom}, together with \eqref{weak.mom.M}, is the weak formulation of \eqref{eq.mom.2}. Also, in \eqref{weak.mom} we define $(\nabla\bm{w})_{ij}\equiv \pa_{x_j}w_i$.
    \end{rem}
    
    \subsection{Existence of weak solutions} \label{ss.ex.weak.sol}
	The goal of this subsection is to prove the following theorem:
	\begin{theo}[Global existence of weak solutions]\label{thm.ex}
		Let $T>0$ be arbitrary, and let the assumptions \eqref{hp.init}--\eqref{hp.eta} hold. Then, a weak solution $(\vel,\mu,\phi,\Phi_v) : Q_T\to \R^3\times\R^n\times\R^n\times\R$ to \eqref{eq.v}--\eqref{eq.Phiv}, \eqref{bc}, \eqref{ic} exists in the sense of Def.~\ref{def.weaksol}.
	\end{theo}	
\begin{rem} 
Standard results in the theory of Navier--Stokes equations~\cite{Tem01} yield the existence of a pressure function $p\in \mathcal{D}'(Q_T)$ such that the momentum balance \eqref{eq.mom.2} holds in the distributional sense.
\end{rem}
    
	The remainder of this paper will be devoted to proving Thm.~\ref{thm.ex}.
	We apply a time semi-discretization and regularization, and solve the resulting stationary elliptic problem via a fix-point argument.

    We introduce the following notation.
    Let $E$ be a bounded measure space and $\mathcal{X}$ a Banach space. Given a sequence $(u_n)_{n\in\N}$ of measurable functions $E\to\mathcal{X}$ and exponents $p>1$, $q<\infty$, we write that $u_n$ is bounded/convergent in $L^{p-}(E; \mathcal{X})$ if $u_n$ is bounded/convergent in $L^{p-\epsilon}(E,\mathcal{X})$ for {\em every} $\epsilon>0$ small enough, while we write that $u_n$ is bounded/convergent in $L^{q+}(E,\mathcal{X})$ if $u_n$ is bounded/convergent in $L^{q+\epsilon}(E; \mathcal{X})$ for {\em some} $\epsilon>0$.

    \medskip
	
	\paragraph{\bf{Approximate problem.}}\label{ss.dt}
	Let $\tau>0$ be an arbitrary (small) time step and $\delta>0$ be a regularization parameter. We denote by $\vec\phi^{k-1}$, $\mathfrak{q}_i^{k-1}$, etc.~the solution at time $k-1$ (which is assumed to be known). Furthermore, recall the definition $\mathfrak{q}_i =\mu_i-\mu_n$, $i=1,\ldots,n-1$.	We split $f^{loc}$ in~\eqref{f} into a convex (w.r.t.~$\vec\phi$) part $g$ and its remainder $h$:
	\begin{align}
		\label{fL.split}
		f^{loc}(\vec\phi,\Phi_v) =& g(\vec\phi) + h(\vec\phi,\Phi_v),\\ 
		g(\vec\phi) =& 
		\sum_{j=1}^n\frac{\phi_j}{N_j}\log\phi_j , \notag\\
		\label{fL.rest}
		h(\vec\phi,\Phi_v) =&
		(1-\omega(\Phi_v))\sum_{i,j=1}^n\chi_{ij}\phi_i\phi_j
        - \omega(\Phi_v)
		\sum_{j=1}^n\frac{\phi_j}{N_j}\log{\phi_j^\mathrm{sat}}.
        \end{align}
	Each part is treated differently in the approximate scheme.	We consider the following time-discrete, regularized version of 
    \eqref{weak.mom}--\eqref{weak.mui}, \eqref{eq.Phiv}, \eqref{bc}, \eqref{ic}: 
	\begin{align}
		\label{app.v}
        &\frac{1}{\tau}\int_\Omega(\rho^k\vel^k - \rho^{k-1}\vel^{k-1})\cdot\bm{w}\dx\\ \nonumber
		&\phantom{=}+\int_\Omega\left(
        2\eta(\vec\phi^k)D\vel^k:D\bm{w}
		+\delta |D\vel^k|^{r-2}D\vel^k:D\bm{w}
		+\delta|\vel^k|^{r-2}\vel^k\cdot\bm{w}
        \right)\dx
		\\ \nonumber
		&=\int_\Omega\Big(
		\rho^k\vel^k\otimes\vel^k 
		+\sum_{j=1}^{n}\beta_j \nabla\phi_j^k\otimes\nabla\phi_j^k
		+\beta_v \nabla\Phi_v^k\otimes\nabla\Phi_v^k
		\Big):D\bm{w}\dx\\ \nonumber
		&\phantom{=}-\int_\Omega\Big( 
		\vel^k\otimes \sum_{i,j=1}^{n-1}(\rhob_i -\rhob_n)(\gamma_{ij}(\vec\phi^k,\Phi_v^k)+\delta\delta_{ij})
		\nabla\mathfrak{q}_j^k\Big):\nabla\bm{w}\dx \quad\forall\bm{w}\in W^{1,r}_\sigma(\Omega),\\ 
		\label{app.phii}
		&\tau\mathfrak{q}_i^k - \Div\sum_{j=1}^{n-1}
		\big(\gamma_{ij}(\vec\phi^k,\Phi_v^k)+\delta\delta_{ij})\nabla\mathfrak{q}_j^k
        +\frac{\phi_i^k-\phi_i^{k-1}}{\tau}
		+\vel^k\cdot\nabla\phi_i^k = 0\\ \nonumber
		&\mbox{in }\Omega,~~i=1,\ldots,n-1,\\ 
		\label{app.mui}
		&-\Delta\sum_{j=1}^{n-1}(\beta_i\delta_{ij}+\beta_n)\phi_j^k
		+\pa_{\phi_i}g(\vec\phi^k)
		+\widetilde{D}_i
		h(\vec\phi^{k},\vec\phi^{k-1},\Phi_v^{k})
		=\mathfrak{q}_i^k\\ \nonumber
		&\mbox{in }\Omega,~~i=1,\ldots,n-1,\\ 
		\label{app.Phiv}
		&-\beta_v\gamma_v(\vec\phi^k)\Delta\Phi_v^k
		+\frac{1}{\tau}(\Phi_v^k - \Phi_v^{k-1}) + 
		\vel^k\cdot\nabla\Phi_v^k
		+\delta (\Phi_v^k)^3\\ \nonumber
		&=-\gamma_v(\vec\phi^k)\frac{f^{loc}(\vec\phi^{k-1},\Phi_v^k)-f^{loc}(\vec\phi^{k-1},\Phi_v^{k-1})}{\Phi_v^k-\Phi_v^{k-1}}
		\qquad\mbox{in }\Omega,
	\end{align}
	where the exponent $r$ in the velocity regularization satisfies $r>4$, and we define
	\begin{align*}
		&\widetilde{D}_i
		h(\vec\phi^{k},\vec\phi^{k-1},\Phi_v^{k}) \\
		\nonumber
		&\equiv  
		\frac{h(\phi_1^{k-1},\ldots,\phi_{i-1}^{k-1},\phi_i^{k},\phi_{i+1}^{k},\ldots,\phi_{n-1}^{k},\Phi_v^{k}) - 
			h(\phi_1^{k-1},\ldots,\phi_{i-1}^{k-1},\phi_i^{k-1},\phi_{i+1}^{k},\ldots,\phi_{n-1}^{k},\Phi_v^{k})}{
			\phi_i^{k}-\phi_i^{k-1}}\\ \nonumber
		&\mbox{for }i=1,\ldots,n-1.
	\end{align*}
	Note that the splitting \eqref{fL.split}--\eqref{fL.rest} is necessary since the terms in $\pa_{\phi_i}g$ are singular in $\vec\phi$ and $g$ is a strictly monotone function of $\vec\phi$. Therefore, $g$ needs to be treated differently from $h$, which is bounded and can be regarded as a mere perturbation. Also, in the above-written expression, we adopt the convention that the difference quotients equal the corresponding partial derivatives when the denominator vanishes.
	
	Define the approximate solution spaces $Y$ and $Y_0$ as follows:
	\begin{align}
		Y &\equiv 
		W^{1,r}_\sigma(\Omega,\R^3)\times H^2(\Omega,\R^{n-1})\times H^2(\Omega), \notag\\
		\label{Y0}
		Y_0 &\equiv \Big\{(\vel,\vec\phi,\Phi_v)\in Y~:~\vec\phi(x)\in\dom \quad\mbox{a.e.~}x\in\Omega\Big\}.
	\end{align}
    We consider regularized initial data $(\vel^{\rm in,\delta},\vec{\phi}^{\rm in,\delta},\Phi_v^{\rm in,\delta})\in Y_0$ in place of the ones defined in \eqref{hp.init}, \eqref{hp.init.2}. Such regularized initial data can be constructed in a standard way.
    
	We derive the following proposition: 
	\begin{prop}[Existence of approximate solutions]\label{prop.ex.app}
		Let $k\geq 1$ be arbitrary.
		Given\\ $(\vel^{k-1},\vec\phi^{k-1},\Phi_v^{k-1})\in Y_0$,
		there exist functions
		$(\vel^{k},\vec\phi^{k},\Phi_v^{k})\in Y_0$, $\vec{\mathfrak{q}}^{k}\in H^2(\Omega;\R^{n-1})$ satisfying \eqref{app.v}--\eqref{app.Phiv}.
	\end{prop}
    The proof of Prop.~\ref{prop.ex.app} is standard and can be found in Appendix~\ref{appendix.proof.prop}. Notice that we implicitly denoted $(\vel^{0},\vec\phi^{0},\Phi_v^{0})\equiv (\vel^{\rm in,\delta},\vec{\phi}^{\rm in,\delta},\Phi_v^{\rm in,\delta})$.

   We continue with the derivation of $\tau$-uniform estimates for the solution $(\vel^k,\vec{\mathfrak{q}}^k,\vec\phi^k,\Phi_v^k)$ to the approximate problem \eqref{app.v}--\eqref{app.Phiv}
	in preparation for the limit $\tau\to 0$. \medskip 
	
	\paragraph{\bf Approximate energy balance.}
    As in the derivation of the energy equation~\eqref{en.bal}, we proceed with a few caveats. Instead of deriving an identity, we derive an inequality in this paragraph. The only differences arise from the difference quotients in the time variable. The following relations hold via elementary convexity arguments:
	\begin{align*}
		(\rho^k\vel^k - \rho^{k-1}\vel^{k-1})\cdot\vel^k - (\rho^k-\rho^{k-1})\frac{|\vel^k|^2}{2}
        &\geq
		\frac{1}{2}\rho^k|\vel^k|^2 - \frac{1}{2}\rho^{k-1}|\vel^{k-1}|^2,\\
		\sum_{i=1}^{n-1}
		(\phi_i^k-\phi_i^{k-1})\pa_{\phi_i}g(\vec\phi^k,\Phi_v^k)
		&\geq 
        g(\vec\phi^k,\Phi_v^k)-g(\vec\phi^{k-1},\Phi_v^k),
	\end{align*}
	while by construction
	\begin{align*}
		\sum_{i=1}^{n-1}
		(\phi_i^k-\phi_i^{k-1})\widetilde{D}_i h(\vec\phi^k,\vec\phi^{k-1},\Phi_v^k) = 
		h(\vec\phi^k,\Phi_v^k)-h(\vec\phi^{k-1},\Phi_v^k).
	\end{align*}
	The transport term yields additional error terms
	\begin{align*}
		&\sum_{i=1}^{n-1}\left(\pa_{\phi_i}g(\vec\phi^k,\Phi_v^k)
		+\widetilde{D}_i
		h(\vec\phi^{k},\vec\phi^{k-1},\Phi_v^{k})\right)\vel^k\cdot\nabla\phi_i^k\\
		&\qquad +\frac{f^{loc}(\vec\phi^{k-1},\Phi_v^k)-f^{loc}(\vec\phi^{k-1},\Phi_v^{k-1})}{\Phi_v^k-\Phi_v^{k-1}}\vel^k\cdot\nabla\Phi_v^k
        =\vel^k\cdot\nabla f^{loc}(\vec\phi^k,\Phi_v^k) + \mathcal{R}^{k,k-1},
	\end{align*}
	where
	\begin{align}\label{err.app.en}
		\mathcal{R}^{k,k-1} = &
		\sum_{i=1}^{n-1}\left(
		\widetilde{D}_i
		h(\vec\phi^{k},\vec\phi^{k-1},\Phi_v^{k}) - 
		\pa_{\phi_i}h(\vec\phi^k,\Phi_v^k)\right)\vel^k\cdot\nabla\phi_i^k\\ \nonumber
		&+\left(
		\frac{f^{loc}(\vec\phi^{k-1},\Phi_v^k)-f^{loc}(\vec\phi^{k-1},\Phi_v^{k-1})}{\Phi_v^k-\Phi_v^{k-1}} - \pa_{\Phi_v}f^{loc}(\vec\phi^k,\Phi_v^k)
		\right)\vel^k\cdot\nabla\Phi_v^k.
	\end{align}	
	In the end, we obtain the {\em approximate energy balance}
	\begin{align}\label{en.bal.app}
		&\int_\Omega H^k\dx + \tau\int_\Omega\Xi_H^k\dx
		+\tau\int_\Omega\mathcal{R}^{k,k-1}\dx
        -\tau^2\int_\Omega\frac{|\vel^k|^2}{2}\sum_{i=1}^{n-1}(\rhob_i-\rhob_n)\mathfrak{q}_i^k\dx 
        \\ \nonumber
		&+\delta\tau\int_\Omega\left(|\vel^k|^r + |D\vel^k|^r + \frac{3}{4}\beta_v |\nabla (\Phi_v^k)^2|^2
		+(\Phi_v^k)^3\frac{f^{loc}(\vec\phi^{k-1},\Phi_v^k)-f^{loc}(\vec\phi^{k-1},\Phi_v^{k-1})}{\Phi_v^k-\Phi_v^{k-1}}
		\right)\dx\\ \nonumber
		&\leq \int_\Omega H^{k-1}\dx
	\end{align}
    and
	\begin{align*}
	&\Xi_H^k\equiv 
	2\eta(\vec\phi^k)|D\vel^k|^2 + \sum_{i,j=1}^{n-1}(\gamma_{ij}(\vec\phi^k,\Phi_v^k)+\delta\delta_{ij})\nabla\mathfrak{q}_i^k\cdot\nabla\mathfrak{q}_j^k + \tau |\vec{\mathfrak{q}}^k|^2
	\\ \nonumber
	&~\qquad +\gamma_v(\vec\phi^k)
	\left|\frac{f^{loc}(\vec\phi^{k-1},\Phi_v^k)-f^{loc}(\vec\phi^{k-1},\Phi_v^{k-1})}{\Phi_v^k-\Phi_v^{k-1}} - \beta_v\Delta\Phi_v^k\right|^2.
	\end{align*}
    From eq.~\eqref{en.bal.app}, we derive several $\tau$-uniform estimates for the approximate solution.
	
	For the time-discrete Allen--Cahn equation \eqref{app.Phiv}, a maximum principle is not available, and we cannot prove uniform $L^\infty$ bounds for $\Phi_v^k$. For this reason, we need to complement \eqref{en.bal.app} with the following estimate, obtained by testing \eqref{app.Phiv} against $\Phi_v^k$:
	\begin{align}\label{app.est.Phiv}
		\delta\tau\|\Phi_v^k\|_{L^4(\Omega)}^4 +&\frac{1}{2}\|\Phi_v^k\|_{L^2(\Omega)}^2\\ \nonumber
		\leq &\frac{1}{2}\|\Phi_v^{k-1}\|_{L^2(\Omega)}^2
		+\frac{\tau}{4}\beta_v^2
		\|\gamma_v(\vec{\phi}^k)^{1/2}\Delta\Phi_v^k\|_{L^2(\Omega)}^2 
		+ C\tau\left(\|\Phi_v^k\|_{L^2(\Omega)}^2 + 1\right).
	\end{align}
	Inequality~\eqref{app.est.Phiv} is needed to control the term involving $(\Phi_v^k)^3$ in~\eqref{en.bal.app} stemming from the regularization.
	
	We sum \eqref{en.bal.app} and \eqref{app.est.Phiv}, and estimate
	\begin{align*}
		\left|\frac{f^{loc}(\vec\phi^{k-1},\Phi_v^k)-f^{loc}(\vec\phi^{k-1},\Phi_v^{k-1})}{\Phi_v^k-\Phi_v^{k-1}} - \beta_v\Delta\Phi_v^k\right|^2\geq 
		\frac{1}{2}\beta_v^2|\Delta\Phi_v^k|^2 - C
	\end{align*}
	to obtain
	\begin{align*}
		&\int_\Omega H^k \dx + \frac{1}{2}\|\Phi_v^k\|_{L^2(\Omega)}^2 
		+c_1\tau\left(
		\|\Delta\Phi_v^k\|_{L^2(\Omega)}^2 
		+\|D\vel^k\|_{L^2(\Omega)}^2
		\right) + c_1\tau^2\|\vec{\mathfrak{q}}^k\|_{L^2(\Omega)}^2\\
		&\qquad + c_1\delta\tau\left(
		\|\nabla\vec{\mathfrak{q}}\|_{L^2(\Omega)}^2
		+\|\vel^k\|_{W^{1,r}(\Omega)}^r + \|\nabla(\Phi_v^k)^2\|_{L^2(\Omega)}^2 
		+ \|\Phi_v^k\|_{L^4(\Omega)}^4
		\right)\\
        &\qquad +\tau^2\int_\Omega\frac{|\vel^k|^2}{2}\sum_{i=1}^{n-1}(\rhob_i-\rhob_n)\mathfrak{q}_i^k\dx \\
		&\leq \int_\Omega H^{k-1}\dx 
		+\frac{1}{2}\|\Phi_v^{k-1}\|_{L^2(\Omega)}^2 
		+\tau\|\mathcal{R}^{k,k-1}\|_{L^1(\Omega)}
		+ C\tau\left(
		\|\Phi_v^{k}\|_{L^2(\Omega)}^2 
		+ 1\right) .
	\end{align*}
    We estimate now via Young's inequality, \eqref{tau.est.q.0} and the fact that {$r>4$}
    \begin{align*}
    \tau^2\left|\int_\Omega\frac{|\vel^k|^2}{2}\sum_{i=1}^{n-1}(\rhob_i-\rhob_n)\mathfrak{q}_i^k\dx \right|\leq &
    \frac{c_1}{2}\tau\delta\|\vel^k\|_{L^r(\Omega)}^r
    +\tau^3\|\mathfrak{q}_i^k\|_{L^2(\Omega)}^2 + C(\delta)\tau\\
    \leq & 
    \frac{c_1}{2}\tau\delta\|\vel^k\|_{L^r(\Omega)}^r + C(\delta)\tau.
    \end{align*}
	From \eqref{err.app.en}, \eqref{E}, and (again) the fact that {$r>4$}, we infer via Young's inequality with a generic parameter $\epsilon>0$
	\begin{align*}
		\|\mathcal{R}^{k,k-1}\|_{L^1(\Omega)}
		&\leq \frac{\delta}{2}\|\vel^k\|_{L^r(\Omega)}^r 
		+ \|\nabla\vec\phi^k\|_{L^2(\Omega)}^2 
		+ \|\nabla\Phi_v^k\|_{L^2(\Omega)}^2 \\
		\nonumber
		&\quad + \|\vec{\phi}^k-\vec{\phi}^{k-1}\|_{L^4(\Omega)}^4
		+ \frac{\epsilon\delta}{8}\|\Phi_v^k-\Phi_v^{k-1}\|_{L^4(\Omega)}^4
		+ C(\delta,\epsilon)\\
		\nonumber
		&\leq \int_\Omega H^k\dx + \frac{\delta}{2}\|\vel^k\|_{L^r(\Omega)}^r 
		+ \epsilon\delta\|\Phi_v^k\|_{L^4(\Omega)}^4
        + \epsilon\delta\|\Phi_v^{k-1}\|_{L^4(\Omega)}^4
		+ C(\delta,\epsilon)
	\end{align*}
	and, by choosing $\epsilon>0$ small enough, we deduce
	\begin{align}\label{app.est.en}
		&\int_\Omega H^k \dx + \frac{1}{2}\|\Phi_v^k\|_{L^2(\Omega)}^2 
		+c_0\tau\left(
		\|\Delta\Phi_v^k\|_{L^2(\Omega)}^2 
		+\|D\vel^k\|_{L^2(\Omega)}^2
		\right) + c_1\tau^2\|\vec{\mathfrak{q}}^k\|_{L^2(\Omega)}^2 \\ \nonumber
		&\qquad + c_0\delta\tau\left(
		\|\nabla\vec{\mathfrak{q}}\|_{L^2(\Omega)}^2
		+\|\vel^k\|_{W^{1,r}(\Omega)}^r + \|\nabla(\Phi_v^k)^2\|_{L^2(\Omega)}^2 
		+ \|\Phi_v^k\|_{L^4(\Omega)}^4
		\right)\\ \nonumber
		&\leq \int_\Omega H^{k-1} \dx + \frac{1}{2}\|\Phi_v^{k-1}\|_{L^2(\Omega)}^2 
		+C\tau\left(
		\int_\Omega H^k \dx + \frac{1}{2}\|\Phi_v^k\|_{L^2(\Omega)}^2 \right)\\ \nonumber
		&\qquad + \frac{c_0\delta\tau}{4}\|\Phi_v^{k-1}\|_{L^4(\Omega)}^4
		+ C(\delta)\tau .
	\end{align}
	\paragraph{\bf Estimates uniform in $\tau$.}
	For every quantity $\{u^k\,:\,k=0,\ldots,m\}\subset \mathcal{X}$ (where $\mathcal{X}$ is any vector space on $\R$) defined on the discrete time grid $\{t_k\,:\,k=0,\ldots,m\}$ with $t_k=k\tau$, $m=T/\tau$, we define a function $u^{(\tau)}:[0,T]\to \mathcal{X}$ such that
	\begin{align*}
		u^{(\tau)}(t) \equiv u^0\chf{\{0\}}(t) 
		+ \sum_{k=1}^m u^{k}\chf{(t_{k-1},t_k]}(t)\quad t\in [0,T],\\
		\hat{u}^{(\tau)}(t)\equiv u^0\chf{\{0\}}(t)
		+ \sum_{k=1}^m\left(
		u^{k-1}\frac{t_k-t}{\tau} + 
		u^{k}\frac{t-t_{k-1}}{\tau}
		\right) \chf{(t_{k-1},t_k]}(t)\quad t\in [0,T],\\
		u^{(\tau,-)}(t) \equiv u^{(\tau)}(t-\tau)\quad t\in [\tau,T].
	\end{align*}
	With this new notation, \eqref{app.v}--\eqref{app.Phiv} become
	\begin{align}
		\label{tau.v}
		&\int_0^T\int_\Omega\pa_t\widehat{\rho\vel}^{(\tau)}\cdot\bm{w}\dx\dt\\ \nonumber
		&\quad +\int_0^T\int_\Omega\left( 2\eta(\vec\phi^{(\tau)})D\vel^{(\tau)}:D\bm{w}
		+\delta|D\vel^{(\tau)}|^{r-2}D\vel^{(\tau)}:D\bm{w}
		+\delta|\vel^{(\tau)}|^{r-2}\vel^{(\tau)}\cdot\bm{w}
        \right)\dx\dt \\ \nonumber
		&=\int_0^T\int_\Omega\Big(
		\rho^{(\tau)}\vel^{(\tau)}\otimes\vel^{(\tau)} 
		+\sum_{j=1}^{n}\beta_j \nabla\phi_j^{(\tau)}\otimes\nabla\phi_j^{(\tau)}
		+\beta_v \nabla\Phi_v^{(\tau)}\otimes\nabla\Phi_v^{(\tau)}
		\Big):D\bm{w}\dx\dt\\ \nonumber
		&\quad -\int_0^T\int_\Omega\Big( 
		\vel^{(\tau)}\otimes \sum_{i,j=1}^{n-1}(\rhob_i -\rhob_n)(\gamma_{ij}(\vec\phi^{(\tau)},\Phi_v^{(\tau)})
		+\delta\delta_{ij})
		\nabla\mathfrak{q}_j^{(\tau)}\Big):\nabla\bm{w}\dx\dt\\ \nonumber
        &\forall\bm{w}\in L^r(0,T; W^{1,r}_\sigma(\Omega)),\\
		& \tau\mathfrak{q}_i^{(\tau)} - \Div \sum_{j=1}^{n-1}
		\big(\gamma_{ij}(\vec\phi^{(\tau)},\Phi_v^{(\tau)})
		+\delta\delta_{ij})\nabla\mathfrak{q}_j^{(\tau)} =
		-\pa_t\hat{\phi}_i^{(\tau)}
		-\vel^{(\tau)}\cdot\nabla\phi_i^{(\tau)} \notag \\ \nonumber
		&\mbox{in }Q_T,~~i=1,\ldots,n-1,\\
		\label{tau.mui}
		&-\Delta\sum_{j=1}^{n-1}(\beta_i\delta_{ij}+\beta_n)\phi_j^{(\tau)}
		+\pa_{\phi_i}g(\vec\phi^{(\tau)})
		+\widetilde{D}_i
		h(\vec\phi^{(\tau)},\vec\phi^{(\tau,-)},\Phi_v^{(\tau)})
		=\mathfrak{q}_i^{(\tau)}\\ \nonumber
		&\mbox{in }Q_T,~~i=1,\ldots,n-1,\\ 
		\label{tau.Phiv}
		&-\beta_v\gamma_v(\vec\phi^{(\tau)})\Delta\Phi_v^{(\tau)}
		+\pa_t\hat{\Phi}_v^{(\tau)}
		+ \vel^{(\tau)}\cdot\nabla\Phi_v^{(\tau)}
		+\delta (\Phi_v^{(\tau)})^3\\ \nonumber
		&=-\gamma_v(\vec\phi^{(\tau)})\frac{f^{loc}(\vec\phi^{(\tau,-)},\Phi_v^{(\tau)})-f^{loc}(\vec\phi^{(\tau,-)},\Phi_v^{(\tau,-)})}{\Phi_v^{(\tau)}-\Phi_v^{(\tau,-)}}
		\quad\mbox{in }Q_T.
	\end{align}
	Summing \eqref{app.est.en} over $k=1,\ldots,m$, using \eqref{E}, and employing a discrete Gronwall inequality leads to (note that $\rho^{(\tau)}$ is uniformly positive and bounded)
	\begin{align*}%
		&\sup_{t\in [0,T]}\left(
		\|\vel^{(\tau)}\|_{L^2(\Omega)}^2 
		+ \|\nabla\vec\phi^{(\tau)}\|_{L^2(\Omega)}^2 
		+ \|\Phi_v^{(\tau)}(t)\|_{H^1(\Omega)}^2
		\right) +\tau\int_0^T\|\vec{\mathfrak{q}}^{(\tau)}\|_{L^2(\Omega)}^2\dt \\ \nonumber
		&+\int_0^T\left(
		\|\Delta\Phi_v^{(\tau)}\|_{L^2(\Omega)}^2 + 
		\|D\vel^{(\tau)}\|_{L^2(\Omega)}^2 + \|\nabla\vec{\mathfrak{q}}^{(\tau)}\|_{L^2(\Omega)}^2
		+\|\vel^{(\tau)}\|_{W^{1,r}(\Omega)}^r
		+\|(\Phi_v^{(\tau)})^2\|_{H^1(\Omega)}^2
		\right)\dt\\ \nonumber
		&\leq C(\delta,T).
	\end{align*}
	The constant $C(\delta,T)$ on the right-hand side of the above inequality also depends on the initial data.
	From the above inequality and the uniform $L^\infty$-bounds for $\vec\phi^{(\tau)}$, we immediately obtain the following bounds:
	\begin{align}\label{tau.est.v}
		\|\vel^{(\tau)}\|_{L^\infty(0,T; L^2(\Omega))}
		+\|D\vel^{(\tau)}\|_{L^2(0,T; L^2(\Omega))}
		+\|\vel^{(\tau)}\|_{L^r(0,T; W^{1,r}(\Omega))}\leq & C(\delta),\\
		\label{tau.est.phii}
		\|\phi_i^{(\tau)}\|_{L^\infty(0,T; L^\infty(\Omega))}
		+\|\phi_i^{(\tau)}\|_{L^\infty(0,T; H^1(\Omega))} \leq & C(\delta),\\
		\|\Phi_v^{(\tau)}\|_{L^\infty(0,T; H^1(\Omega))}
		+\|\Phi_v^{(\tau)}\|_{L^2(0,T; H^2(\Omega))}
		+\| (\Phi_v^{(\tau)})^2 \|_{L^2(0,T; H^1(\Omega))} \leq & C(\delta), \notag\\
		\label{tau.est.q.0}
        \sqrt{\tau}\|\mathfrak{q}\|_{L^2(0,T; L^2(\Omega))} + 
		\|\nabla\mathfrak{q}\|_{L^2(0,T; L^2(\Omega))} \leq & C(\delta).
	\end{align}
    Inequality~\eqref{fp.est.3a} yields the following bound 
    $$
    \int_0^T\left(  
    \int_\Omega(\phi_i^{(\tau)}-\xi_i^{(\tau)})\log\frac{\phi_i^{(\tau)}}{\xi_i^{(\tau)}}\dx    
    \right)^2\dt\leq C(\delta),\quad i=1,\ldots,n,
    $$
    with $\xi_i^{(\tau)}(t) = \frac{1}{|\Omega|}\int_\Omega\phi_i^{(\tau)}(t)\dx$. Integrating \eqref{app.phii} yields 
    $$
    \int_\Omega\phi_i^{k}\dx = \int_\Omega\phi_i^{k-1}\dx
    -\tau^2\int_\Omega\mathfrak{q}_i^k\dx,\quad i=1,\ldots,n-1,
    $$
    thus summing the above identity in $k$ yields
    \begin{align*}
    \int_\Omega\phi_i^{(\tau)}(t)\dx &= \int_\Omega\phi_i^{{\rm in},\delta}\dx - \tau\int_0^t\int_\Omega\mathfrak{q}_i^{(\tau)}\dx\dt ,\quad i=1,...,n-1,\\
    \int_\Omega\phi_n^{(\tau)}(t)\dx &= \int_\Omega\phi_n^{{\rm in},\delta}\dx + \tau\int_0^t\int_\Omega\sum_{j=1}^{n-1}\mathfrak{q}_j^{(\tau)}\dx\dt.
    \end{align*}
    From the above identity, the a.e.~positivity of $\phi_i^{{\rm in},\delta}$ and \eqref{tau.est.q.0} we infer that $\xi_{\rm min}>0$, $\tau_0(\delta)>0$ exists such that $\xi_i^{(\tau)}(t) \geq \xi_{\rm min}>0$ a.e.~$t\in [0,T]$, $i=1,\ldots,n$, provided that $\tau<\tau_0(\delta)$. Therefore from \eqref{fp.est.3a}, \eqref{fp.est.3b}, and the above estimates we obtain 
	\begin{align}\label{tau.est.phii.2}
		\|\phi_i^{(\tau)}\|_{L^2(0,T; H^2(\Omega))}
		+\|\log\phi_i^{(\tau)}\|_{L^2(0,T; L^1(\Omega))} \leq C(\delta),\quad i=1,\ldots,n .
	\end{align}
	Thus, one easily obtains a $\tau$-uniform estimate for $\int_\Omega\mathfrak{q}_i\dx$ by integrating \eqref{app.mui} over $\Omega$. Consequently,
	\begin{align}\label{tau.est.q}
		\|\mathfrak{q}^{(\tau)}\|_{L^2(0,T; H^1(\Omega))} \leq  C(\delta).
	\end{align}
    We will make use of the following elementary inequality. Given $\vec\phi\in\mathring{\dom}$, there exists $\hat{j}\in\{1,...,n\}$ such that $\phi_{\hat{j}}\geq\frac{1}{n}$. Therefore,
    $$
    \frac{1}{N_i} | \log\phi_i | \leq \left|\frac{1}{N_i}\log\phi_i - 
    \frac{1}{N_{\hat{j}}}\log\phi_{\hat{j}}\right| + \frac{1}{N_{\hat{j}}}|\log\phi_{\hat{j}}|\leq 
    \sum_{k=1}^n \left|\frac{1}{N_i}\log\phi_i-\frac{1}{N_k}\log\phi_k\right| + \frac{n}{\min_{j}N_j}
    $$
    for $i=1,...,n$.
	From \eqref{tau.mui}, \eqref{tau.est.phii.2}, \eqref{tau.est.q}, and the above inequality we infer
	\begin{align}\label{tau.est.logphii}
		\|\log\phi_i^{(\tau)}\|_{L^2(0,T; L^2(\Omega))}\leq C(\delta),\quad i=1,\ldots,n,
	\end{align}
	while from \eqref{tau.est.phii.2} and \eqref{tau.est.v} it follows that
	\begin{align*}
		\|\rho^{(\tau)}\vel^{(\tau)}\|_{L^2(0,T; H^1(\Omega))}\leq C (\delta).
	\end{align*}
	Lastly, using \eqref{tau.v}--\eqref{tau.Phiv}, we deduce
	\begin{align*}
		\|\pa_t(\widehat{\rho\vel}^{(\tau)})\|_{L^{\frac{r}{r-1}}(0,T; (W^{1,r}_\sigma(\Omega))')} + 
		\|\pa_t\hat{\phi}_i^{(\tau)}\|_{L^2(0,T; (H^1(\Omega))')}
		+\|\pa_t\hat{\Phi}_v^{(\tau)}\|_{L^{\frac{4}{3}}(0,T; L^{\frac{4}{3}}(\Omega))}\leq C(\delta).
	\end{align*}
	\paragraph{\bf Passing to the limit $\tau\to 0$.}
	From the $\tau$-uniform estimates and Aubin--Lions lemma we deduce the following convergence relations:
	\begin{align*}
		\rho^{(\tau)}\vel^{(\tau)}\to\rho\vel,
		\quad \phi_i^{(\tau)}\to\phi_i,\quad 
		\Phi_v^{(\tau)}\to\Phi_v
		\quad\mbox{strongly in }L^2(0,T; L^2(\Omega)).
	\end{align*}
	This estimate can be improved by employing the uniform estimates from the previous subsection, as well as the following Gagliardo--Nirenberg--Sobolev inequality:
	\begin{align*}
		\|\nabla\Phi_v^{(\tau)}\|_{L^{16/5}(\Omega)}\leq 
		C\|\Phi_v^{(\tau)}\|_{L^8(\Omega)}^{1/2}
		\|\Phi_v^{(\tau)}\|_{H^2(\Omega)}^{1/2}.
	\end{align*}
	We obtain
	\begin{align}\label{tau.cnv.v}
		\vel^{(\tau)}\to\vel\quad &\mbox{strongly in }L^{r-}(0,T; L^\infty(\Omega)),\\
		\phi_i^{(\tau)}\to\phi_i\quad &\mbox{strongly in }
        L^p(Q_T)\cap L^{2-}(0,T; W^{1,6-}(\Omega))\mbox{ for every }p<\infty, \notag\\
		\Phi_v^{(\tau)}\to\Phi_v\quad
		&\mbox{strongly in }L^{8-}(Q_T)\cap L^{16/5-}(0,T; W^{1,16/5-}(\Omega)), \notag\\
		\vel^{(\tau)}\rightharpoonup\vel\quad
		&\mbox{weakly in }L^r(0,T; W^{1,r}(\Omega)), \notag\\
		\phi_i^{(\tau)}\rightharpoonup\phi_i\quad
		&\mbox{weakly in }L^2(0,T; H^2(\Omega)), \notag\\
		\Phi_v^{(\tau)}\rightharpoonup\Phi_v\quad
		&\mbox{weakly in }L^2(0,T; H^2(\Omega))\cap L^8(Q_T), \notag\\
        \tau\vec{\mathfrak{q}}\to 0\quad & \mbox{strongly in }L^2(0,T; L^2(\Omega)). \notag
	\end{align}
	Finally, \eqref{tau.est.logphii} implies $\vec\phi\in\dom$ and
	\begin{align}\label{tau.cnv.logphii}
		\log\phi_i^{(\tau)}\to \log\phi_i\quad\mbox{strongly in }L^{2-}(Q_T).
	\end{align}
	At this point, it is straightforward to take the limit $\tau\to 0$ in \eqref{tau.v}--\eqref{tau.Phiv} and obtain the following time-continuous regularized equations:
	\begin{align}
	\label{delta.v}
&\int_0^T\int_\Omega\pa_t(\rho^{(\delta)}\vel^{(\delta)})\cdot\bm{w}\dx\dt\\ \nonumber
	& +\int_0^T\int_\Omega\left(
    2\eta(\vec\phi^{(\delta)})D\vel^{(\delta)}:D\bm{w}
	+\delta|D\vel^{(\delta)}|^{r-2}D\vel^{(\delta)}:D\bm{w}
	+\delta|\vel^{(\delta)}|^{r-2}\vel^{(\delta)}\cdot\bm{w}
    \right)\dx\dt	\\ \nonumber
	&=\int_0^T\int_\Omega\Big(
	\rho^{(\delta)}\vel^{(\delta)}\otimes\vel^{(\delta)} 
	+\sum_{j=1}^{n}\beta_j \nabla\phi_j^{(\delta)}\otimes\nabla\phi_j^{(\delta)}
	+\beta_v \nabla\Phi_v^{(\delta)}\otimes\nabla\Phi_v^{(\delta)}
	\Big):D\bm{w}\dx\dt\\ \nonumber
	& -\int_0^T\int_\Omega\Big( 
	\vel^{(\delta)}\otimes \sum_{i,j=1}^{n-1}(\rhob_i -\rhob_n)(\gamma_{ij}(\vec\phi^{(\delta)},\Phi_v^{(\delta)})+\delta\delta_{ij})
	\nabla\mathfrak{q}_j^{(\delta)}\Big):\nabla\bm{w}\dx\dt\\ \nonumber
    &\forall\bm{w}\in L^r(0,T; W_\sigma^{1,r}(\Omega)),\\
	\label{delta.phii}
	&\pa_t{\phi}_i^{(\delta)}
	+\vel^{(\delta)}\cdot\nabla\phi_i^{(\delta)} =
	\Div \sum_{j=1}^{n-1}
	\big(\gamma_{ij}(\vec\phi^{(\delta)},\Phi_v^{(\delta)})
	+\delta\delta_{ij})\nabla\mathfrak{q}_j^{(\delta)}
	\\ \nonumber
	&\mbox{in }Q_T,~~i=1,\ldots,n-1,\\
	&-\Delta\sum_{j=1}^{n-1}(\beta_i\delta_{ij}+\beta_n)\phi_j^{(\delta)}
	+\pa_{\phi_i}f^{loc}(\vec\phi^{(\delta)},\Phi_v^{(\delta)})
	=\mathfrak{q}_i^{(\delta)} \notag\\ \nonumber
	&\mbox{in }Q_T,~~i=1,\ldots,n-1,\\ 
	\label{delta.Phiv}
	&-\beta_v\gamma_v(\vec\phi^{(\delta)})\Delta\Phi_v^{(\delta)}
	+\pa_t\Phi_v^{(\delta)}
	+ \vel^{(\delta)}\cdot\nabla\Phi_v^{(\delta)}
	+\delta (\Phi_v^{(\delta)})^3\\ \nonumber
	&=-\gamma_v(\vec\phi^{(\delta)})
	\pa_{\Phi_v}f^{loc}(\vec\phi^{(\delta)},\Phi_v^{(\delta)})
	\quad\mbox{in }Q_T.
	\end{align}
    In order to take the limit $\tau\to 0$ in \eqref{en.bal.app}, let us consider the term $\mathcal{R}^{k,k-1}$. From \eqref{err.app.en} we obtain
    \begin{align*}
    \tau &\sum_{k=1}^m \int_\Omega\mathcal{R}^{k,k-1}\dx \\ 
    =& \int_0^T\int_\Omega
    \sum_{i=1}^{n-1}\left(\widetilde{D}_i
	h(\vec\phi^{(\tau)},\vec\phi^{(\tau,-)},\Phi_v^{(\tau)}) - 
	\pa_{\phi_i}h(\vec\phi^{(\tau)},\Phi_v^{(\tau)})\right)\vel^{(\tau)}\cdot\nabla\phi_i^{(\tau)}
     \dx\dt \\ \nonumber
	&+\int_0^T\int_\Omega\left(
		\frac{f^{loc}(\vec\phi^{(\tau,-)},\Phi_v^{(\tau)})-f^{loc}(\vec\phi^{(\tau,-)},\Phi_v^{(\tau,-)})}{\Phi_v^{(\tau)}-\Phi_v^{(\tau,-)}} - \pa_{\Phi_v}f^{loc}(\vec\phi^{(\tau)},\Phi_v^{(\tau)})
		\right)\vel^{(\tau)}\cdot\nabla\Phi_v^{(\tau)}\dx\dt .
    \end{align*}
    Given the definitions of $\widetilde{D}h$, $h$, $f^{loc}$, it is easy to see that 
    \begin{align*}
    \widetilde{D}_i
	h(\vec\phi^{(\tau)},\vec\phi^{(\tau,-)},\Phi_v^{(\tau)}) - 
	\pa_{\phi_i}h(\vec\phi^{(\tau)},\Phi_v^{(\tau)})\to 0\quad\mbox{a.e.~in }Q_T,\\
	\frac{f^{loc}(\vec\phi^{(\tau,-)},\Phi_v^{(\tau)})-f^{loc}(\vec\phi^{(\tau,-)},\Phi_v^{(\tau,-)})}{\Phi_v^{(\tau)}-\Phi_v^{(\tau,-)}} - \pa_{\Phi_v}f^{loc}(\vec\phi^{(\tau)},\Phi_v^{(\tau)})\to 0\quad\mbox{a.e.~in }Q_T,
    \end{align*}
    and the functions are bounded in $L^\infty(Q_T)$, while $\vel^{(\tau)}\cdot\nabla\Phi_v^{(\tau)}$, $\vel^{(\tau)}\cdot\nabla\phi_i^{(\tau)}$ are bounded in $L^{1+}(Q_T)$ thanks to the previously shown uniform bounds. It follows that $\tau \sum_{k=1}^m \int_\Omega\mathcal{R}^{k,k-1}\dx\to 0$ as $\tau\to 0$.
    This property and the convergence relations
    \eqref{tau.cnv.v}--\eqref{tau.cnv.logphii} allow us to sum \eqref{en.bal.app} in $k=1,...,\overline{k}$ with $\overline{k}=\lfloor t/\tau\rfloor$ and $t\in (0,T)$ arbitrary, take the limit $\tau\to 0$ in the resulting inequality and obtain the following {\em regularized energy balance}:
	\begin{align}\label{en.bal.delta}
		&\int_\Omega H^{(\delta)}(t)\dx + \int_0^t\int_\Omega\Xi_H^{(\delta)}\dx\dt'
		\\ \nonumber
		&\qquad +\delta\int_0^t\int_\Omega\left(|\vel^{(\delta)}|^r + |D\vel^{(\delta)}|^r + \frac{3}{4}\beta_v|\nabla (\Phi_v^{(\delta)})^2|^2
        +(\Phi_v^{(\delta)})^3
		\pa_{\Phi_v}f^{loc}(\vec\phi^{(\delta)},\Phi_v^{(\delta)})
		\right)\dx\dt'\\ \nonumber
		&\leq \int_\Omega H^{\mathrm{in}}\dx 
	\end{align}
	with
	\begin{align}
		\label{en.bal.Xi.delta}
		&\Xi_H^{(\delta)}\equiv 
		2\eta(\vec\phi^{(\delta)})|D\vel^{(\delta)}|^2 + \sum_{i, j=1}^{n-1}(\gamma_{ij}(\vec\phi^{(\delta)},\Phi_v^{(\delta)})+\delta\delta_{ij})\nabla\mathfrak{q}_i^{(\delta)}\cdot\nabla\mathfrak{q}_j^{(\delta)} 
		\\ \nonumber
		&~\qquad +\gamma_v(\vec\phi^{(\delta)})
		\left|\pa_{\Phi_v}
		f^{loc}(\vec\phi^{(\delta)},\Phi_v^{(\delta)})
		- \beta_v\Delta\Phi_v^{(\delta)}\right|^2.
	\end{align}
	\medskip
	
	\paragraph{\bf Estimates uniform in $\delta$.}
    Passing to the limit $\tau\to 0$ in \eqref{tau.est.phii}, \eqref{tau.est.phii.2} yields ($\delta$-dependent) bounds for $\phi_i^{(\delta)}$ in $L^\infty(0,T; L^\infty(\Omega))$ and $L^2(0,T; H^2(\Omega))$. Given that the mobility $\gamma_v\in C^1(\dom)$, it follows via a Gagliardo--Nirenberg argument that $|\nabla \gamma_v(\vec\phi^{(\delta)})|\in L^4(0,T; L^4(\Omega))$. 

It is possible to prove that 
    \begin{align}\label{delta.bound.PhiV01}
		0\leq\Phi_v^{(\delta)}\leq 1\quad\mbox{a.e.~in }Q_T.
	\end{align}
    If the estimate $|\nabla \gamma_v(\vec\phi^{(\delta)})|\in L^{\frac{2p}{p-3}}(0,T; L^p(\Omega))$ held for some $p>3$, then \eqref{delta.bound.PhiV01} would follow from a straightforward Stampacchia truncation argument applied to \eqref{delta.Phiv} and the identity
    $-\beta_v\gamma_v(\vec\phi^{(\delta)})\Delta\Phi_v^{(\delta)} = 
    -\Div(\beta_v\gamma_v(\vec\phi^{(\delta)})\nabla\Phi_v^{(\delta)})
    +\beta_v\nabla \gamma_v(\vec\phi^{(\delta)})\cdot \nabla\Phi_v^{(\delta)}$. However, the aforementioned regularity estimate for $\nabla\gamma_v(\vec\phi^{(\delta)})$ does not hold in general. A workaround is the following argument. First, one considers \eqref{delta.Phiv} with $\gamma_v(\vec\phi^{(\delta)})$ replaced by a suitable approximation $\gamma_v^{(\ell)}\in L^{\frac{2p}{p-3}}(0,T; W^{1,p}(\Omega))$, $\ell\in\N$, such that $\gamma_v^{(\ell)}\to \gamma_v(\vec\phi^{(\delta)})$ strongly in $L^4(0,T; W^{1,4}(\Omega))$: 
    \begin{align}\label{AC.ell}
    &\pa_t\Phi_v^{(\delta,\ell)}
    -\beta_v\gamma_v^{(\ell)}\Delta\Phi_v^{(\delta,\ell)}
	+ \vel^{(\delta)}\cdot\nabla\Phi_v^{(\delta,\ell)}
	+\delta (\Phi_v^{(\delta,\ell)})^3
    =-\gamma_v^{(\ell)}\pa_{\Phi_v}f^{loc}(\vec\phi^{(\delta)},\Phi_v^{(\delta,\ell)})\mbox{ in }Q_T.    
    \end{align}
    The corresponding solution $\Phi_v^{(\delta,\ell)}$ to \eqref{AC.ell} would fulfill \eqref{delta.bound.PhiV01}. Standard arguments would then yield the convergence of $\Phi_v^{(\delta,\ell)}$ towards a solution to \eqref{delta.Phiv}. Finally, standard energy methods applied to \eqref{delta.Phiv} would yield the uniqueness of solutions to \eqref{delta.Phiv}. As a consequence, each solution $\Phi_v^{(\delta)}$ to \eqref{delta.Phiv} would satisfy \eqref{delta.bound.PhiV01}. As the full argument would be rather technical but quite standard, we omit it and just assume \eqref{delta.bound.PhiV01} holds.
	We define the chemical potentials $\mu_i^{(\delta)}$ as 
	\begin{align*}
		\mu_i^{(\delta)} = \pa_{\phi_i}f^{loc}(\phi_1^{(\delta)},\ldots,\phi_n^{(\delta)},\Phi_v^{(\delta)}) - \beta_i\Delta\phi_i^{(\delta)}\quad i=1,\ldots,n,
	\end{align*}
	so that $\mathfrak{q}_i^{(\delta)} = \mu_i^{(\delta)}-\mu_n^{(\delta)}$ for $i=1,\ldots,n-1$. Also, the bounds \eqref{est.LinfL2}--\eqref{est.namu} hold, namely,
	\begin{align}\label{est.LinfL2.delta}
		\|\nabla\phi_i^{(\delta)}\|_{L^\infty(0,T; L^2(\Omega))} + 
		\|\nabla\Phi_v^{(\delta)}\|_{L^\infty(0,T; L^2(\Omega))} +
		\|\vel^{(\delta)}\|_{L^\infty(0,T; L^2(\Omega))} &\leq C,\\
		\label{est.L2L2.delta}
		\|D\vel^{(\delta)}\|_{L^2(0,T; L^2(\Omega))} 
		+ \|\Delta\Phi_v^{(\delta)}\|_{L^2(0,T; L^2(\Omega))} &\leq C,\\
		\int_0^T\int_\Omega \sum_{k,j=1}^{n-1}\gamma_{kj}^{(\delta)}\nabla(\mu_k^{(\delta)}-\mu_n^{(\delta)})\cdot\nabla(\mu_j^{(\delta)}-\mu_n^{(\delta)})\dx\dt &\leq C, \notag
	\end{align}	
	for $i=1,\ldots,n-1$, with the additional bound
	\begin{align*}
		\delta\int_0^T\int_\Omega|\nabla\mathfrak{q}_i^{(\delta)}|^2\dx\dt &\leq C.
	\end{align*}
	Therefore, we can apply Lemma~\ref{lem.est.namu} and obtain
	\begin{align}\label{delta.bound.sqrtphii}
		\|\nabla\sqrt{\phi_i^{(\delta)}}\|_{L^2(0,T; L^2(\Omega))} &\leq C,\\
		\label{delta.bound.Deltaphii}
		\|D^2{\phi_i^{(\delta)}}\|_{L^2(0,T; L^2(\Omega))} &\leq C,\\
		\int_0^T\int_\Omega\sum_{i=1}^{n}\phi_i^{(\delta)}\left|
		\nabla\mu_i^{(\delta)} - \sum_{j=1}^{n}\alpha_j^{(\delta)}\nabla\mu_j^{(\delta)}
		\right|^2\dx\dt &\leq C,
		\qquad\alpha_i^{(\delta)} \equiv \frac{D_i\phi_i^{(\delta)}}{\sum_{k=1}^n D_k\phi_k^{(\delta)}}. \notag%
	\end{align}
	In particular, $(\phi_i^{(\delta)})^\epsilon\mu_i^{(\delta)}$ is bounded in $L^2(0,T; L^2(\Omega))$ for every $\epsilon>0$. Furthermore, it holds that
	\begin{align*}
		\|J^{V,i,(\delta)}\|_{L^2(0,T; L^2(\Omega))} + 
		\|J^{L,i,(\delta)}\|_{L^2(0,T; L^2(\Omega))}\ &\leq C,
	\end{align*}
	where $J^{V,i,(\delta)}$ and $J^{L,i,(\delta)}$ are the vapor and liquid components of the regularized diffusion flux, i.e.
	\begin{align*}
		J^{V,i,(\delta)} &= -\omega(\Phi_v^{(\delta)})\sum_{j=1}^{n-1}
		\gamma_{ij}^V(\vec\phi^{(\delta)}) \nabla(\mu_j^{(\delta)}-\mu_n^{(\delta)})
        =-\omega(\Phi_v^{(\delta)})
		\sum_{j=1}^n\gamma_{ij}^V(\vec\phi^{(\delta)})
		\nabla\mu_j^{(\delta)},\\
		J^{L,i,(\delta)} &= -(1-\omega(\Phi_v^{(\delta)}))\sum_{j=1}^{n-1}
		\gamma_{ij}^L
		(\vec\phi^{(\delta)})
		\nabla(\mu_j^{(\delta)}-\mu_n^{(\delta)})
        = -(1-\omega(\Phi_v^{(\delta)}))\sum_{j=1}^{n}
		\gamma_{ij}^L
		(\vec\phi^{(\delta)})\nabla\mu_j^{(\delta)}.
	\end{align*}
	The total regularized diffusion flux $\bJ^{(\delta)} = J^{L,(\delta)} + J^{V,(\delta)}$ satisfies
	\begin{align}\label{bound.delta.J}
		\|\bJ^{(\delta)}\|_{L^2(Q_T)}\leq C.
	\end{align}
	For arbitrary $\Psi_1,\ldots,\Psi_n\in L^2(0,T; H^1_0(\Omega,\R^3))$, integration by parts yields
	\begin{align}
		\label{Ji.delta}
		\sum_{i=1}^n &\int_0^T\int_\Omega\Psi_i\cdot \bJ_i^{(\delta)}\dx\dt =
		\sum_{i,j=1}^n\int_0^T \int_\Omega\gamma_{ij}^{(\delta)}\mu_j^{(\delta)}\,\Div\Psi_i\dx\dt\\ \nonumber
		&+\sum_{i,j=1}^n\int_0^T \int_\Omega\Psi_i\cdot
		\left(
		\sum_{\ell=1}^n\frac{(\phi_\ell^{(\delta)})^{\epsilon_0}}{1-\epsilon_0}\,\pa_{\phi_\ell}\gamma_{ij}^{(\delta)}\nabla (\phi_\ell^{(\delta)})^{1-\epsilon_0}
		+\pa_{\Phi_v}\gamma_{ij}^{(\delta)}\nabla\Phi_v^{(\delta)}
		\right)\mu_j^{(\delta)}\dx\dt
	\end{align}
	with $\gamma_{ij}^{(\delta)} = \gamma_{ij}^{(\delta)}(\phi_1^{(\delta)},\ldots,\phi_n^{(\delta)},\Phi_v^{(\delta)})$ and $\epsilon_0\in (0,1/2)$ as in \eqref{hp.mob.1}.
	We point out that \eqref{est.Lp} holds, namely
	\begin{align}
		\label{est.Lp.delta}
		\|\nabla\Phi_v^{(\delta)}\|_{L^{10/3}(Q_T)}+
		\|\nabla\phi_i^{(\delta)}\|_{L^{10/3}(Q_T)}+
		\|\vel^{(\delta)}\|_{L^{10/3}(Q_T)}\leq C,\quad i=1,...,n.
	\end{align}
	From \eqref{delta.phii} and \eqref{bound.delta.J} we deduce
	\begin{align}\label{est.phiit.delta}
		\|\pa_t\phi_i^{(\delta)}\|_{L^2(0,T; (H^1(\Omega))')}\leq C.
	\end{align}
	From \eqref{est.LinfL2.delta}, \eqref{est.Lp.delta}, \eqref{est.phiit.delta}, and Aubin--Lions lemma we infer (up to subsequences)
	\begin{align*}
		\phi_i^{(\delta)}\to\phi_i\quad\mbox{strongly in }L^{10/3}(0,T; C^0(\overline{\Omega}))\cap C^0([0,T]; L^2(\Omega)).
	\end{align*}
	Given the uniform $L^\infty$-bounds for $\phi$ we deduce
	\begin{align}
		\label{phii.cnv.delta}
		&\phi_i^{(\delta)}\to\phi_i\quad\mbox{strongly in }L^{p}(0,T; L^\infty(\Omega))\cap L^{\infty}(0,T; L^p(\Omega))\cap C^0([0,T]; L^2(\Omega))\\ \nonumber
        &\forall p<\infty,\quad i=1,\ldots,n.
	\end{align}
	With a similar argument we deduce the strong convergence of $\Phi_v^{(\delta)}$,
	\begin{align}
		\label{PhiV.cnv.delta}
		&\Phi_v^{(\delta)}\to\Phi_v\quad\mbox{strongly in }L^{p}(0,T; L^\infty(\Omega))\cap L^{\infty}(0,T; L^p(\Omega))\cap C^0([0,T]; L^2(\Omega))\\ \nonumber
        &\forall p<\infty.
	\end{align}
	The above convergence relations, bounds~\eqref{est.L2L2.delta} and~\eqref{delta.bound.Deltaphii},
	and the Gagliardo--Nirenberg--Sobolev inequality lead to
	\begin{align}
		\label{delta.cnv.gradphi}
		\nabla\phi_i^{(\delta)}\to \nabla\phi_i,\quad 
		\nabla\Phi_v^{(\delta)}\to \nabla\Phi_v\quad 
		\mbox{strongly in }L^{\frac{10}{3}-}(Q_T).
	\end{align}
	From \eqref{delta.v} and the previous estimates we infer
	\begin{align}
	\label{delta.rhovt}
	\|\pa_t(\rho^{(\delta)}\vel^{(\delta)})\|_{L^{s/(s-1)}(0,T; (W^{1,s}_\sigma(\Omega))')}\leq C,\quad\mbox{for some }s>1.
	\end{align}
	On the other hand, from \eqref{est.Lp.delta} we obtain
	\begin{align}\label{delta.narhov}
	\|\rho^{(\delta)}\vel^{(\delta)}\|_{L^{10/3}(Q_T)} + 
	\|\nabla(\rho^{(\delta)}\vel^{(\delta)})\|_{L^{5/3}(Q_T)}\leq C.
	\end{align}
	Therefore, via Aubin--Lions lemma, we obtain 
	the strong convergence of $\rho^{(\delta)}\vel^{(\delta)}$ in $L^{5/3}(Q_T)$. Furthermore, $\rho^{(\delta)}$ is strongly convergent, uniformly positive, and bounded, hence \eqref{est.Lp.delta} implies
	\begin{align*}
	\vel^{(\delta)}\to\vel \quad\mbox{strongly in }L^{\frac{10}{3}-}(Q_T).
	\end{align*}
	We now investigate the convergence of the chemical potentials. Let us define
	\begin{align*}
	\mu_i = \begin{cases}
\pa_{\phi_i}f^{loc}(\phi_1,\ldots,\phi_n,\Phi_v) - \beta_i\Delta\phi_i & \mbox{on }Q_T\cap\{\phi_i>0\},\\
0 & \mbox{on }Q_T\cap\{\phi_i=0\},	
	\end{cases}
	\quad i=1,\ldots,n.
	\end{align*}
	Given the previous convergence results and \eqref{delta.bound.Deltaphii} we deduce
	\begin{align}\label{delta.cnv.mu}
	(\phi_i^{(\delta)})^\epsilon\mu_i^{(\delta)}\rightharpoonup (\phi_i)^\epsilon\mu_i\quad\mbox{weakly in }L^2(Q_T),~~\forall\epsilon>0.
	\end{align}
    In order to pass to the limit $\delta\to 0$ in \eqref{Ji.delta} we rewrite 
    \begin{align*}
	(\phi_\ell^{(\delta)})^{\epsilon_0}\,\pa_{\phi_\ell}\gamma_{ij}^{(\delta)}\nabla (\phi_\ell^{(\delta)})^{1-\epsilon_0} \mu_j^{(\delta)} = 
    (\phi_\ell^{(\delta)})^{\epsilon_0}\,\pa_{\phi_\ell}\gamma_{ij}^{(\delta)}
    (\phi_j^{(\delta)})^{-\epsilon_0/2}\cdot 
    \nabla (\phi_\ell^{(\delta)})^{1-\epsilon_0} \cdot 
    (\phi_j^{(\delta)})^{\epsilon_0/2}\mu_j^{(\delta)} 
    \end{align*}
    From Assumption \eqref{hp.mob.1} and the strong convergence relations \eqref{phii.cnv.delta}, \eqref{PhiV.cnv.delta} we deduce that the first factor in the above decomposition is a.e.~convergent in $Q_T$ and bounded in $L^\infty(Q_T)$. The second factor $\nabla (\phi_\ell^{(\delta)})^{1-\epsilon_0}$ is strongly convergent in $L^2(Q_T)$ thanks to \eqref{delta.cnv.gradphi}, \eqref{delta.bound.sqrtphii} and Lemma~\ref{lem.aux.1}. Finally, the third factor $(\phi_j^{(\delta)})^{\epsilon_0/2}\mu_j^{(\delta)}$ is weakly convergent in $L^2(Q_T)$ thanks to \eqref{delta.cnv.mu}.
    We can therefore take the limit $\delta\to 0$ in \eqref{Ji.delta}
	and obtain \eqref{Ji}.
	
	Using the previously proven convergence relations, we are also able to pass to the limit $\delta\to 0$ in 
	\eqref{delta.v}--\eqref{delta.Phiv} and deduce
	\eqref{weak.mom}--\eqref{weak.phii},
	\eqref{eq.Phiv}. 
	
	Finally, rewriting the term
	\begin{align*}
	\sum_{i,j=1}^{n-1}\gamma_{ij}(\vec\phi^{(\delta)},\Phi_v^{(\delta)})\nabla\mathfrak{q}_i^{(\delta)}\cdot\nabla\mathfrak{q}_j^{(\delta)}  = &
	\sum_{i,j=1}^{n}\gamma_{ij}(\vec\phi^{(\delta)},\Phi_v^{(\delta)})\nabla\mu_i^{(\delta)}\cdot
	\nabla\mu_j^{(\delta)} \\ 
    =& \mathcal{Q}(\vec\phi^{(\delta)},\Phi_v^{(\delta)})[\bJ_1^{(\delta)},\ldots,\bJ_n^{(\delta)}]
	\end{align*}
    with $\mathcal{Q}$ as in \eqref{Q} 
    and employing the convergence relations allows us to pass to the limit $\delta\to 0$ in \eqref{en.bal.delta}, \eqref{en.bal.Xi.delta},
	yielding \eqref{weak.en.bal}, \eqref{weak.en.bal.Xi} for $t_1=0$. However, it is possible to obtain \eqref{weak.en.bal} for a.e.~$0\leq t_1<t_2\leq T$ by summing \eqref{en.bal.app} in $k=k_1,...,k_2$, with $k_1 = \lfloor t_1/\tau\rfloor$, $k_2=\lfloor t_2/\tau\rfloor$, then take the limit $\tau\to 0$ in the resulting inequality, yielding an analogue inequality to \eqref{en.bal.delta} but with $t$ replaced by $t_2$, $\int_0^t$ replaced by $\int_{t_1}^{t_2}$, and $H^{\rm in}$ replaced by $H^{(\delta)}(t_2)$. This is possible since the convergence relations proved till \eqref{en.bal.delta} yield the strong $L^1(\Omega)$ convergence of $H^{(\tau)}(t)$ for a.e.~$t\in (0,T)$. Then we take the limit $\delta\to 0$ in the new inequality to obtain \eqref{weak.en.bal} for a.e.~$0\leq t_1<t_2\leq T$, again by using the strong $L^1(\Omega)$ convergence of $H^{(\delta)}(t)$ for a.e.~$t\in (0,T)$.

    We now show that the initial conditions are satisfied in the sense of \eqref{ic.weak}. From \eqref{delta.narhov}, \eqref{delta.rhovt}, \eqref{phii.cnv.delta}, \eqref{PhiV.cnv.delta} we infer in a standard way that 
    $\vel\in C^0_{\rm weak}([0,T]; L^2(\Omega))$, $\phi_i, \Phi_v\in C^0([0,T]; L^2(\Omega))$, and $\phi_i(0)=\phi_i^{\mathrm{in}}$, $\Phi_v(0)=\Phi_v^{\mathrm{in}}$, $\vel(0)=\vel^{\mathrm{in}}$.
    From the energy balance \eqref{weak.en.bal} we deduce that 
    \begin{align*}
    \limsup_{t\to 0^+} &\int_\Omega\left(
    \rho(t) |\vel(t)|^2 + \sum_{j=1}^n\beta_j|\nabla\phi_j(t)|^2 + 
    \beta_v|\nabla\Phi_v(t)|^2\right)\dx\\
    \leq &
    \int_\Omega\left(
    \rho(0) |\vel(0)|^2 + \sum_{j=1}^n\beta_j|\nabla\phi_j(0)|^2 + 
    \beta_v|\nabla\Phi_v(0)|^2\right)\dx,
    \end{align*}
    implying the strong $L^2(\Omega)$-convergence of $\sqrt{\rho(t)}\vel(t)$, $\nabla\phi_i(t)$, $\nabla\Phi_v(t)$ towards $\sqrt{\rho(0)}\vel(0)$, $\nabla\phi_i(0)$, $\nabla\Phi_v(0)$ as $t\to 0$ (respectively). 
    Since $(\sqrt{\rho(t)}-\sqrt{\rho^{\mathrm{in}}})\vel^{\mathrm{in}}\to 0$ strongly in $L^2(\Omega)$ via dominated convergence and $\sqrt{\rho(t)}\geq (\min_{i=1,...,n}\rhob_i)^{1/2}$, we infer that 
    \begin{align*}
    \|\vel(t)-\vel^{\mathrm{in}}\|_{L^2(\Omega)} &\leq C\|\sqrt{\rho(t)}(\vel(t)-\vel^{\mathrm{in}})\|_{L^2(\Omega)}\\ 
    &\leq 
    C\|\sqrt{\rho(t)}\vel(t)-\sqrt{\rho^{\mathrm{in}}}\vel^{\mathrm{in}}\|_{L^2(\Omega)} + 
    C\|(\sqrt{\rho(t)}-\sqrt{\rho^{\mathrm{in}}})\vel^{\mathrm{in}}\|_{L^2(\Omega)}\to 0
    \end{align*}
    as $t\to 0^+$. We conclude that \eqref{ic.weak} holds.
	This completes the proof of Thm.~\ref{thm.ex}.

	\appendix
	\section{Proof of Prop.~\ref{prop.ex.app}}
    \label{appendix.proof.prop}
		The proof is based on a fix point strategy. 
    Define the approximate solution spaces $X$, $X_0$, $\mathcal{B}$ as follows:
	\begin{align*}
		X &\equiv 
		W^{1,r}_\sigma(\Omega,\R^3)\times H^{2}(\Omega,\R^{n-1})\times H^2(\Omega,\R^{n-1})\times H^2(\Omega),\\
		X_0 &\equiv \Big\{(\vel,\vec{\mathfrak{q}},\vec\phi,\Phi_v)\in X~:~\vec\phi(x)\in\dom \quad\mbox{a.e.~}x\in\Omega\Big\},\\
		\mathcal{B} &\equiv L^\infty_\sigma(\Omega,\R^3)\times W^{1,4}(\Omega,\R^{n-1})\times W^{1,4}(\Omega,\R^{n-1})\times W^{1,4}(\Omega).
		\end{align*}
	Note that $X\hookrightarrow \mathcal{B}$ compactly (Sobolev's embedding theorem). Next, define the projection mapping
		\begin{align*}
			\mathcal{P}_{\dom}:\R^{n-1}\to\dom,\qquad 
			\mathcal{P}_{\dom}(\vec\phi)_i := 
			\frac{(\phi_i)_+}{ \max\{1,\sum_{j=1}^{n-1}(\phi_j)_+\} }\quad (i=1,\ldots,n-1).
		\end{align*}
		For $(\vel^*,\vec{\mathfrak{q}}^*,\vec\phi^*,\Phi_v^*)\in\mathcal{B}$, 
        $(\vel^{k-1},\vec\phi^{k-1},\Phi_v^{k-1})\in Y_0$ (recall \eqref{Y0}),
        consider the PDE system
		\begin{align}
			\label{fp.v}
			&\frac{1}{\tau}\int_\Omega(\rho^*\vel - \rho^{k-1}\vel^{k-1})\cdot\bm{w}\dx\\ \nonumber
			&\quad +\int_\Omega\left(
            2\eta(\mathcal{P}_{\dom}(\vec\phi^*))D\vel:D\bm{w} 
            +\delta|D\vel|^{r-2}D\vel:D\bm{w}
            +\delta|\vel|^{r-2}\vel\cdot\bm{w}\right)\dx
			\\ \nonumber
			&=\int_\Omega\Big(
			\rho^*\vel^*\otimes\vel^* 
			+\sum_{j=1}^{n}\beta_j \nabla\phi_j^*\otimes\nabla\phi_j^*
			+\beta_v \nabla\Phi_v^*\otimes\nabla\Phi_v^*
			\Big):D\bm{w}\dx\\ \nonumber
			&\quad -\int_\Omega\Big( 
			\vel^*\otimes \sum_{i,j=1}^{n-1}(\rhob_i -\rhob_n)(\gamma_{ij}(\mathcal{P}_{\dom}(\vec\phi^*),\Phi_v^*)+\delta\delta_{ij})
			\nabla\mathfrak{q}_j^*\Big):\nabla\bm{w}\dx\\ \nonumber
            & \forall\bm{w}\in W^{1,r}_\sigma(\Omega),\\
			\label{fp.phii}
			&\tau\mathfrak{q}_i
            -\Div \sum_{j=1}^{n-1}
			\big(\gamma_{ij}(\mathcal{P}_{\dom}(\vec\phi^*),\Phi_v^*)+\delta\delta_{ij})\nabla\mathfrak{q}_j =
			-\frac{1}{\tau}(\phi_i^*-\phi_i^{k-1})
			-\vel^*\cdot\nabla\phi_i^* \\ \nonumber
			&\mbox{in }\Omega,~~i=1,\ldots,n-1,\\
			\label{fp.mui}
			&-\Delta\sum_{j=1}^{n-1}(\beta_i\delta_{ij}+\beta_n)\phi_j 
			+\pa_{\phi_i}g(\vec\phi)
			+\widetilde{D}_i
			h(\vec\phi^{*},\vec\phi^{k-1},\Phi_v^{*})
			=\mathfrak{q}_i^*\\ \nonumber
			&\mbox{in }\Omega,~~i=1,\ldots,n-1,\\ 
			\label{fp.Phiv}
			&-\beta_v\gamma_v(\mathcal{P}_{\dom}(\vec\phi^*))\Delta\Phi_v
			+\frac{1}{\tau}(\Phi_v - \Phi_v^{k-1})  
			+\delta\Phi_v^3\\ \nonumber
			&=-\vel^*\cdot\nabla\Phi_v^*
			-\gamma_v(\mathcal{P}_{\dom}(\vec\phi^*))
			\frac{f^{loc}(\vec\phi^{k-1},\Phi_v^*)-f^{loc}(\vec\phi^{k-1},\Phi_v^{k-1})}{\Phi_v^*-\Phi_v^{k-1}}
			\qquad\mbox{in }\Omega,
		\end{align}
		where $\rho^* = \rhob_n + \sum_{i=1}^{n-1}(\rhob_i-\rhob_n)\mathcal{P}_{\dom}(\vec\phi^*)_i$.		
		We couple the system with the boundary conditions~\eqref{bc}.
		The remainder of the proof is divided into two steps.\medskip\\
		{\em Step 1: Well-posedness of the fix point map.}
		We show that \eqref{fp.v}--\eqref{fp.Phiv} has a unique solution $(\vel,\vec{\mathfrak{q}},\vec\phi,\Phi_v)\in\mathcal{B}$. 
		
		Eq.~\eqref{fp.v} is a stationary incompressible Stokes equation with an added monotone term $-\delta\Div(|D\vel|^{r-2}D\vel)+\delta|\vel|^{r-2}\vel$ and a right-hand side in $(W^{1,r}_\sigma(\Omega))'$; the existence of a unique weak solution $\vel\in W^{1,r}_\sigma(\Omega)$ is a classical matter~\cite{Tem01}.
		
		Eq.~\eqref{fp.phii} admits a unique weak solution $\vec{\mathfrak{q}}\in H^1(\Omega)$. However, given that $\vel^*\in L^{\infty}_\sigma(\Omega)$, it holds that $\vel^*\cdot\nabla\phi_i^* = \Div(\phi_i^*\vel^*)$ and $\phi_i^*\vel^*\in L^\infty(\Omega)$. Hence standard elliptic regularity theory yields $\vec{\mathfrak{q}}\in W^{1,p}(\Omega)$ for every $p<\infty$. With this new information in mind, we write \eqref{fp.phii} in non-divergence form, thus obtaining (up to multiplication by a bounded, symmetric matrix) the Poisson equation for $\mathfrak{q}_i$ ($i=1,\ldots,n-1$) and derive once again via elliptic regularity results that (at least) $\vec{\mathfrak{q}}\in H^2(\Omega)$ (in fact, $\vec{\mathfrak{q}}\in W^{2,p}(\Omega)$ but the $H^2$-bound is sufficient in this setting).
		
		Eq.~\eqref{fp.Phiv} is also solved in a standard way to obtain a solution $\Phi_v\in H^2(\Omega)$.
		
		Next, we wish to solve \eqref{fp.mui}. The main difficulty in solving the equation lies in the singularity in the term 
		\begin{align}\label{fp.monotone.repr}
			\pa_{\phi_i}g (\vec\phi,\Phi_v^{*})
			=& \pa_{\phi_i}\left[\sum_{j=1}^{n-1}N_j^{-1}\psi(\phi_j) 
			+N_n^{-1}\psi\left(1-\sum_{j=1}^{n-1}\phi_j\right)
			\right],
		\end{align}
		where $\psi(s) = s\log s$ for $s\geq 0$.
		To regularize the singular terms in \eqref{fp.monotone.repr}, define for $0<\epsilon<\frac{1}{e}$
		\begin{align*}
			\psi_\epsilon(s) = 
			\begin{cases}
				s\log s & s\geq \epsilon,\\
				\epsilon\log\epsilon + (\log\epsilon + 1)(s-\epsilon) + \frac{1}{2\epsilon}(s-\epsilon)^2 & s<\epsilon.
			\end{cases}
		\end{align*}
		The function $\psi_\epsilon\in C^2(\R)$ is a regularization of $\psi$. We replace $\psi$ by $\psi_\epsilon$ in \eqref{fp.monotone.repr} to obtain a regularized version of \eqref{fp.mui}:
		\begin{align}
			\label{fp.monot.eps}
			-\Delta &\sum_{j=1}^{n-1}B_{ij}\phi_j^{(\epsilon)}
			+N_i^{-1}\psi_\epsilon'(\phi_i^{(\epsilon)}) - N_n^{-1}
			\psi_\epsilon'\left(1-\sum_{j=1}^{n-1}\phi_j^{(\epsilon)}\right)
			=\mathfrak{q}_i^*
			-\widetilde{D}_i h(\vec\phi^{*},\vec\phi^{k-1},\Phi_v^{*}),
		\end{align}
		where we define $B_{ij}=\beta_i\delta_{ij}+\beta_n$.
		
		Since $\psi_\epsilon : \R\to\R$ is convex and $B$ is symmetric, positive definite, and constant, eq.~\eqref{fp.monot.eps} has the structure $Au=b$, where $A: H^1(\Omega)\to (H^1(\Omega))'$ is a strictly monotone, coercive, and hemicontinuous operator. This can be solved classically~\cite[Thm.~26A]{Zei2B}, yielding a weak solution $\vec\phi^{(\epsilon)}\in H^1(\Omega,\R^{n-1})$. Standard elliptic regularity results yield $\vec\phi^{(\epsilon)}\in H^2(\Omega,\R^{n-1})$.
		The convexity of $\psi_\epsilon$ yields an
		$\epsilon$-uniform $H^1$-bound for $\vec\phi^{(\epsilon)}$ (by testing \eqref{fp.monot.eps} with $\phi_i^{(\epsilon)}$). Furthermore, we can derive an $\epsilon$-uniform $H^2$-bound for $\vec\phi^{(\epsilon)}$ by testing \eqref{fp.monot.eps} with $-\Delta\phi_i^{(\epsilon)}$ and integrating the monotone term by parts. More precisely, the contribution of the latter reads
		\begin{align*}
			-\int_\Omega & N_i^{-1}\psi_\epsilon'(\phi_i^{(\epsilon)})\Delta\phi_i^{(\epsilon)}
            \dx
			=\int_\Omega N_i^{-1}\psi_\epsilon''(\phi_i^{(\epsilon)})|\nabla\phi_i^{(\epsilon)}|^2\dx
			\geq -C.
		\end{align*}
		Thus, we obtain
		\begin{align*}
			\|\vec\phi^{(\epsilon)}\|_{H^2(\Omega)}\leq C.
		\end{align*}
		In particular, the monotonicity of $\psi_\epsilon'$ yields, upon testing \eqref{fp.monot.eps} with $\phi_i^{(\epsilon)}-\xi_i$, where $\vec\xi\in\mathring{\dom}$ is an arbitrary constant,
		\begin{align}
			\label{fp.monot.bound.2}
			\int_\Omega(\phi_i^{(\epsilon)}-\xi_i)(\psi_\epsilon'(\phi_i^{(\epsilon)})-\psi_\epsilon'(\xi_i))\dx\leq C,\quad i=1,\ldots,n,
		\end{align}
		which in turn implies
		\begin{align*}
			C\geq &\int_{\Omega\cap\{\phi_i^{(\epsilon)}\leq 0\}}(\phi_i^{(\epsilon)}-\xi_i)\psi_\epsilon'(\phi_i^{(\epsilon)})\dx
			=\int_{\Omega\cap\{\phi_i^{(\epsilon)}\leq 0\}}(\xi_i+|\phi_i^{(\epsilon)}|)(\log(\epsilon^{-1})+\epsilon^{-1} |\phi_i^{(\epsilon)}| ) 	\dx\\ 
			\geq &\frac{1}{\epsilon}\int_{\Omega\cap\{\phi_i^{(\epsilon)}\leq 0\}}|\phi_i^{(\epsilon)}|^2\dx,\quad i=1,\ldots,n.
		\end{align*}
		The compact Sobolev embedding $H^2(\Omega)\hookrightarrow L^\infty({\Omega})$ yields the strong convergence $\vec\phi^{(\epsilon)}\to\vec\phi$ in $L^\infty(\Omega)$. The above relation implies that $\vec\phi\in{\dom}$ a.e.~in $\Omega$. However, Fatou's lemma, applied to \eqref{fp.monot.bound.2}, allows us to infer that
        \begin{align*}
			\int_\Omega(\phi_i-\xi_i)(\log(\phi_i)-\log(\xi_i))\dx\leq C,\quad i=1,\ldots,n,
		\end{align*}  
        and hence, $\vec\phi\in\mathring{\dom}$ a.e.~in $\Omega$. Furthermore, 
        $\vec\phi\in H^2(\Omega)\hookrightarrow C^{0,\alpha}(\overline{\Omega})$. 
        In particular, the sets $\mathcal{S}_k:=\cap_{i=1}^n\{x\in\Omega\,:\,\phi_i(x)>2^{-k}\}$ are open and $\cup_{k\in\N}\mathcal{S}_k = \Omega$ up to a null set. Hence by taking arbitrary test functions $\eta\in C^\infty_c$ supported in some $\mathcal{S}_k$, we can take the limit $\epsilon\to 0$ in \eqref{fp.monot.eps} and obtain \eqref{fp.mui}.
		
		In conclusion, \eqref{fp.v}--\eqref{fp.Phiv} has a unique solution $(\vel,\vec{\mathfrak{q}},\vec\phi,\Phi_v) \in X\subset\mathcal{B}$, and we thus obtain a mapping 
		$$
		\fun{\mathcal{T}}{\mathcal{B}}{\mathcal{B}}{(\vel^*,\vec{\mathfrak{q}}^*,\vec\phi^*,\Phi_v^*)}{(\vel,\vec{\mathfrak{q}},\vec\phi,\Phi_v)}.
		$$
		{\em Step 2: The mapping $\mathcal{T}$ has a fix point.}
		From the arguments in the previous step we deduce that $(\vel,\vec{\mathfrak{q}},\vec\phi,\Phi_v)$ is bounded in $X$ whenever $(\vel^*,\vec{\mathfrak{q}}^*,\vec\phi^*,\Phi_v^*)$ is bounded in $\mathcal{B}$, following from the standard bounds in classical regularity theory for elliptic equations. This means that $\mathcal{T}:\mathcal{B}\to\mathcal{B}$ is a compact operator via standard Sobolev embeddings. It is also standard to show the (sequential) continuity of $\mathcal{T}:\mathcal{B}\to\mathcal{B}$.
		We will employ Schaefer's fix point theorem to show the existence of a fix point for $\mathcal{T}$. To this end, we consider the set
		$$
		\{ (\vel^*,\vec{\mathfrak{q}}^*,\vec\phi^*,\Phi_v^*)\in\mathcal{B}~~:~~
		(\vel^*,\vec{\mathfrak{q}}^*,\vec\phi^*,\Phi_v^*) = \sigma\mathcal{T}(\vel^*,\vec{\mathfrak{q}}^*,\vec\phi^*,\Phi_v^*) \}
		$$
		and show it is bounded in $\mathcal{B}$ uniformly w.r.t.~$\sigma\in [0,1]$. Let $(\vel,\vec{\mathfrak{q}},\vec\phi,\Phi_v)=\mathcal{T}(\vel^*,\vec{\mathfrak{q}}^*,\vec\phi^*,\Phi_v^*)$. Hence, \eqref{fp.v}--\eqref{fp.Phiv} hold with 
		$(\vel^*,\vec{\mathfrak{q}}^*,\vec\phi^*,\Phi_v^*)=\sigma(\vel,\vec{\mathfrak{q}},\vec\phi,\Phi_v)$.
		Notice that $\mathcal{P}_{\dom}(\vec\phi^*)=\mathcal{P}_{\dom}(\sigma\vec\phi) = \sigma\vec\phi$ since $\vec\phi(\Omega)\subset\dom$. 
		
		Given that $\vel^*\in L^\infty_\sigma(\Omega)$, it holds that $\vel^*\cdot\nabla\phi_i^* = \Div(\phi_i^*\vel^*)$ and $\phi_i^*\vel^*\in L^\infty(\Omega)$ (via Sobolev's embedding), hence standard energy methods applied to \eqref{fp.phii} yield the bound
		\begin{align}\label{fp.est.1}
        \|\mathfrak{q}_i\|_{L^2(\Omega)} + 
		\|\nabla\mathfrak{q}_i\|_{L^2(\Omega)}\leq C(\delta,\tau)(1+\|\vel\|_{L^2(\Omega)}).
		\end{align}
		Choosing $\bm{w}=\vel$ in \eqref{fp.v} and employing Korn's inequality leads to
		\begin{align*}
			\|\vel\|_{W^{1,r}(\Omega)}^r
			\leq &C(\delta,\tau)\Bigg(
			1 + \int_\Omega|D\vel||\vel|^2\dx
			+ \int_\Omega\sum_{i=1}^{n-1}|\nabla\phi_i|^2 |D\vel|\dx \\  
			&~\qquad\qquad + \int_\Omega|\nabla\Phi_v|^2 |D\vel|\dx
			+ \int_\Omega|\nabla\vel||\vel||\nabla\vec{\mathfrak{q}}|\dx
			\Bigg),
		\end{align*}
		therefore, via Young's inequality, we infer
		\begin{align*}
			&\|\vel\|_{W^{1,r}(\Omega)}^r
			\leq C(\delta,\tau)\left(
			1 + \sum_{i=1}^{n-1}\int_\Omega |\nabla\phi_i|^{\frac{2r}{r-1}}\dx
			+ \int_\Omega |\nabla\Phi_v|^{\frac{2r}{r-1}}\dx
			+ \int_\Omega |\nabla\vec{\mathfrak{q}}|^{\frac{r}{r-2}}\dx
			\right).
		\end{align*}
		Since $\frac{r}{r-2}<2$ for $r>4$, from \eqref{fp.est.1} we deduce that the integral involving $\nabla\vec{\mathfrak{q}}$ can be absorbed (via Young's inequality) by the left-hand side of the above inequality, yielding
		\begin{align}\label{fp.est.2}
			&\|\vel\|_{W^{1,r}(\Omega)}^r
			\leq C(\delta,\tau)\left(
			1 + \sum_{i=1}^{n-1}\int_\Omega |\nabla\phi_i|^{\frac{2r}{r-1}}\dx
			+ \int_\Omega |\nabla\Phi_v|^{\frac{2r}{r-1}}\dx
			\right).
		\end{align}
		Next, we derive a bound for $\vec\phi$ in $H^1$.
		Let $\xi_i\equiv\frac{1}{|\Omega|}\int_\Omega\phi_i\dx$. We employ $\phi_i - \xi_i$ as a test function in \eqref{fp.mui} to obtain 
		\begin{align}\label{fp.est.3.1}
			\sum_{i=1}^{n-1}\|\nabla\phi_i\|_{L^2(\Omega)}^2
			+\sum_{i=1}^{n}\int_\Omega(\phi_i-\xi_i)\log\frac{\phi_i}{\xi_i}\dx
			\leq C_1 	
			+C_2\sum_{i=1}^{n-1}\int_\Omega\mathfrak{q}_i(\phi_i-\xi_i)\dx.
		\end{align}
		Since $\int_\Omega(\phi_i-\xi_i)\dx = 0$ we can rewrite
		$$
		\int_\Omega\mathfrak{q}_i(\phi_i-\xi_i)\dx = 
		\int_\Omega\left(\mathfrak{q}_i-|\Omega|^{-1}\int_\Omega\mathfrak{q}_i\dx \right)(\phi_i-\xi_i)\dx
		$$
		and, via Cauchy--Schwarz and Poincar\'e--Wirtinger inequalities, we deduce from \eqref{fp.est.3.1} the following bound:
		\begin{align}\label{fp.est.3a}
			\sum_{i=1}^{n-1}\|\nabla\phi_i\|_{L^2(\Omega)}^2
			+\sum_{i=1}^{n}\int_\Omega(\phi_i-\xi_i)\log\frac{\phi_i}{\xi_i}\dx
			\leq C \left(1 +\sum_{i=1}^{n-1}\|\nabla\mathfrak{q}_i\|_{L^2(\Omega)} \right).
		\end{align}
		At this point, to derive a bound for $\vec\phi$ in $H^2$, we employ $\Delta\phi_i$ as a test function in \eqref{fp.mui} and obtain
		\begin{align}\label{fp.est.3b}
			\sum_{i=1}^{n-1}\|\phi_i\|_{H^2(\Omega)}^2\leq C(\delta)\left(1 + \|\nabla\Phi_v\|_{L^2(\Omega)}^2 + \sum_{i=1}^{n-1}\|\nabla\mathfrak{q}_i\|_{L^2(\Omega)}^2
			\right).
		\end{align}
		Summing \eqref{fp.est.3a}, \eqref{fp.est.3b} and employing \eqref{fp.est.1} leads to
		\begin{align}\label{fp.est.3}
			\sum_{i=1}^{n-1} &\|\phi_i\|_{H^2(\Omega)}^2
			+\sum_{i=1}^{n}\int_\Omega(\phi_i-\xi_i)\log\frac{\phi_i}{\xi_i}\dx
			\leq C(\delta,\tau)\left(1 + \|\nabla\Phi_v\|_{L^2(\Omega)}^2 +\|\vel\|_{L^2(\Omega)}^2 
			\right).
		\end{align}
		We derive a $\sigma$-uniform bound for $\Phi_v$ in $W^{1,4}(\Omega)$ (note that one could show even higher regularity).
		Multiplying \eqref{fp.Phiv} by $\Phi_v$ and integrating over $\Omega$ yields
		\begin{align*}
			\frac{1}{4\tau}\int_\Omega\Phi_v^2\dx + \delta\int_\Omega\Phi_v^4\dx\leq 
			\frac{1}{2\tau}\int_\Omega|\Phi_v^{k-1}|^2\dx
			+C_0\tau\int_\Omega(\Delta\Phi_v)^2\dx+C
		\end{align*}
		for some $\tau$-independent constant $C_0>0$. Next,	
		multiplying \eqref{fp.Phiv} by $-\Delta\Phi_v$ and integrating over $\Omega$ yields
		\begin{align*}
			c_0\int_\Omega &(\Delta\Phi_v)^2\dx + \frac{1}{2\tau}\int_\Omega|\nabla\Phi_v|^2\dx 
			+\frac{3}{4}\delta\int_\Omega|\nabla(\Phi_v^2)|^2\dx\\
			\leq & C + \frac{1}{2\tau}\int_\Omega|\nabla\Phi_v^{k-1}|^2\dx 
			+\sigma^2\int_\Omega\Delta\Phi_v(\vel\cdot\nabla)\Phi_v\dx
		\end{align*}
		since \eqref{hp.gammaV} holds and the factor multiplied by $\gamma_v(\vec\phi^*)$ on the right-hand side of \eqref{fp.Phiv} is bounded. 
		Given that $\vel$ is divergence-free, it follows, via integration by parts, that
		\begin{align*}
			\sigma^2\int_\Omega\Delta\Phi_v(\vel\cdot\nabla)\Phi_v\dx = 
			-\sigma^2\int_\Omega\nabla\Phi_v\cdot (D\vel) \nabla\Phi_v\dx,
		\end{align*}
		hence, by summing the previous two inequalities while assuming w.l.o.g.~$\tau\leq\frac{c_0}{2C_0}$, we deduce
		\begin{align}\label{fp.est.Phiv.1}
			\|\Delta\Phi_v\|_{L^2(\Omega)}^2 + \|\Phi_v\|_{H^1(\Omega)}^2 + \|\Phi_v^2\|_{H^1(\Omega)}^2 
			\leq C_1(\delta,\tau)\int_\Omega |D\vel||\nabla\Phi_v|^2\dx + C(\delta,\tau).
		\end{align}
		We now estimate
		\begin{align*}
			C_1(\delta,\tau)
			\int_\Omega |D\vel||\nabla\Phi_v|^2\dx\leq 
			\epsilon\|D\vel\|_{L^r(\Omega)}^r + 
			\epsilon_1\|\nabla\Phi_v\|_{L^{8/3}(\Omega)}^{8/3} + C(\epsilon,\epsilon_1,\delta,\tau),
		\end{align*}
		where we use that {$r>4$} to apply Young's inequality. The Gagliardo--Nirenberg--Sobolev inequality,
		\begin{align}\label{GNS}
			\|\nabla\Phi_v\|_{L^{8/3}(\Omega)}\leq 
			C\|\Phi_v\|_{L^{4}(\Omega)}^{1/2}
			\|\Phi_v\|_{H^2(\Omega)}^{1/2},
		\end{align}
		and Young's inequality allow us to infer
		\begin{align*}
			C_1(\delta,\tau)\int_\Omega |D\vel||\nabla\Phi_v|^2\dx\leq 
			\epsilon\|D\vel\|_{L^r(\Omega)}^r + 
			\epsilon_1\|\Phi_v\|_{L^{4}(\Omega)}^{4} 
			+\epsilon_1\|\Phi_v\|_{H^2(\Omega)}^{2}
			+ C(\epsilon,\epsilon_1,\delta,\tau).
		\end{align*}
		Estimating $\|\Phi_v\|_{L^{4}(\Omega)}^{4} = 
		\|\Phi_v^2\|_{L^{2}(\Omega)}^{2}\leq C\|\Phi_v^2\|_{H^{1}(\Omega)}^{2}$ 
		and choosing $\epsilon_1>0$ small enough allows us to infer from \eqref{fp.est.Phiv.1} the following estimate:
		\begin{align}\label{fp.est.Phiv}
			\|\Delta\Phi_v\|_{L^2(\Omega)}^2 + \|\Phi_v\|_{H^1(\Omega)}^2 + \|\Phi_v^2\|_{H^1(\Omega)}^2 
			\leq C(\delta,\tau,\epsilon) + \epsilon\|D\vel\|_{L^r(\Omega)}^r.
		\end{align}	
		We sum \eqref{fp.est.2}, \eqref{fp.est.3}, \eqref{fp.est.Phiv}, and choose $\epsilon>0$ small enough to obtain
		\begin{align}\label{fp.est.fin}
			&\|\vel\|_{W^{1,r}(\Omega)}^r 
			+\sum_{i=1}^{n-1}
            \left( 
            \|\phi_i\|_{H^2(\Omega)}^2
			+\int_\Omega(\phi_i-\xi_i)\log\frac{\phi_i}{\xi_i}\dx
            \right) \\ \nonumber
			&\qquad +\|\Delta\Phi_v\|_{L^2(\Omega)}^2 + \|\Phi_v\|_{H^1(\Omega)}^2 + \|\Phi_v^2\|_{H^1(\Omega)}^2 
			\\ \nonumber
			&\leq C(\delta,\tau)\left(
			1 + \sum_{i=1}^{n-1}\int_\Omega |\nabla\phi_i|^{\frac{2r}{r-1}}\dx
			+ \int_\Omega |\nabla\Phi_v|^{\frac{2r}{r-1}}\dx
			\right)	.
		\end{align}
		Since {$r>4$}, it follows that $\frac{2r}{r-1}<\frac{8}{3}$. Thus, Young's inequality, the Gagliardo--Nirenberg--Sobolev inequality~\eqref{GNS}, and the uniform $L^\infty$-bounds for $\vec\phi$ allow us to estimate the right-hand side of \eqref{fp.est.fin} in terms of the left-hand side. We conclude that there is a $\sigma$-uniform bound for $(\vel,\vec{\mathfrak{q}},\vec\phi,\Phi_v)$
		(and a fortiori for $(\vel^*,\vec{\mathfrak{q}}^*,\vec\phi^*,\Phi_v^*)=\sigma(\vel,\vec{\mathfrak{q}},\vec\phi,\Phi_v)=\sigma\mathcal{T}(\vel^*,\vec{\mathfrak{q}}^*,\vec\phi^*,\Phi_v^*)$) in $\mathcal{B}$. Therefore, Schaefer's fix point theorem yields the existence of a fix point $(\vel^k,\vec{\mathfrak{q}}^k,\vec\phi^k,\Phi_v^k)\in\mathcal{B}$ for $\mathcal{T}$, i.e.,~functions $(\vel^k,\vec\phi^k,\Phi_v^k)\in Y_0$, $\vec{\mathfrak{q}}^k\in H^2(\Omega;\R^{n-1})$ satisfying \eqref{app.v}--\eqref{app.Phiv}. This completes the proof of the proposition.

    \section{Auxiliary results}\label{app.auxres}
    We prove here some miscellaneous results employed in the analysis.
    \begin{lemma}\label{lem.aux.1}
    Let $(\phi_k)_{k\in\N}\subset L^2(0,T; H^1(\Omega))$ be a sequence such that $0\leq\phi_k\leq 1$ a.e.~in $Q_T$, $\phi_k\to\phi$ strongly in $L^2(0,T; H^1(\Omega))$ and
    $$
    \exists p>2:\quad 
    \int_{Q_T}(|\nabla\phi_k|^p + |\nabla\phi_k^{1/2}|^2)\dx\dt\leq C.
    $$
    Then, for every $\alpha\in (0,1/2)$ it holds that 
    $\phi_k^{\alpha+1/2}\to \phi^{\alpha+1/2}$ strongly in $L^2(0,T; H^1(\Omega))$.
    \end{lemma}
    \begin{proof}
We define the exponent $q\in (2,p)$ such that $\frac{1}{q}=\frac{1-2\alpha}{2}+\frac{2\alpha}{p}$ and notice that $|\nabla\phi_k^{\alpha+1/2}|^q = c|\nabla\phi_k|^{2\alpha q}|\nabla\sqrt{\phi_k}|^{(1-2\alpha)q}$ for some suitable constant $c$ (depending on $\alpha$, $q$). Therefore, H\"older's inequality yields
\begin{align*}
\int_{Q_T}|\nabla\phi_k^{\alpha+1/2}|^q\dx\dt \leq \left(
\int_{Q_T}|\nabla\phi_k|^p\dx\dt \right)^{\frac{2\alpha q}{p}}
\left(
\int_{Q_T}|\nabla\sqrt{\phi_k}|^2\dx\dt \right)^{\frac{(1-2\alpha)q}{2}}\leq C.
\end{align*}
We are now going to show that $\nabla\phi_k^{\alpha+1/2}\to\nabla\phi^{\alpha+1/2}$ strongly in $L^1(Q_T)$.
For arbitrary $\eta>0$, define $Q_{T,\eta}\equiv \{(x,t)\in Q_{T}~:~\phi(x,t)\geq\eta\}$. By Egorov's theorem, for every $\epsilon>0$ there exists $Q_{T,\eta,\epsilon}\subset Q_{T,\eta}$ such that $|Q_{T,\eta}\backslash Q_{T,\eta,\epsilon}|<\epsilon$, $\phi_k\to\phi$, $\nabla\phi_k\to\nabla\phi$ uniformly in $Q_{T,\eta,\epsilon}$, and $\phi_k\geq \eta/2$ a.e.~in $Q_{T,\eta,\epsilon}$. We split
\begin{align*}
\int_{Q_T}|\nabla\phi_k^{\alpha+1/2}-\nabla\phi^{\alpha+1/2}|\dx\dt = \mathcal{I}_{1,k} + \mathcal{I}_{2,k} + \mathcal{I}_{3,k},\\
\mathcal{I}_{1,k} = \int_{Q_{T,\eta,\epsilon}}|\nabla\phi_k^{\alpha+1/2}-\nabla\phi^{\alpha+1/2}|\dx\dt,\\
\mathcal{I}_{2,k} = \int_{Q_{T,\eta}\backslash Q_{T,\eta,\epsilon}}|\nabla\phi_k^{\alpha+1/2}-\nabla\phi^{\alpha+1/2}|\dx\dt,\\
\mathcal{I}_{3,k} = \int_{Q_{T}\backslash Q_{T,\eta}}|\nabla\phi_k^{\alpha+1/2}-\nabla\phi^{\alpha+1/2}|\dx\dt.
\end{align*}
Clearly $\mathcal{I}_{1,k}\to 0$ as $k\to\infty$. Furthermore
\begin{align*}
\mathcal{I}_{2,k}\leq &|Q_{T,\eta}\backslash Q_{T,\eta,\epsilon}|^{1-1/q}\left(
\|\nabla\phi_k^{\alpha+1/2}\|_{L^q(Q_T)} + \|\nabla\phi^{\alpha+1/2}\|_{L^q(Q_T)}\right)\leq C\epsilon^{1-1/q}.
\end{align*}
Finally, the $L^2$ bound for $\nabla\sqrt{\phi_k}$ yields
\begin{align*}
\mathcal{I}_{3,k}\leq &C\int_{Q_{T}\backslash Q_{T,\eta}}\phi_k^{\alpha}|\nabla\sqrt{\phi_k}|\dx\dt + C\int_{Q_{T}\backslash Q_{T,\eta}}\phi^{\alpha}|\nabla\sqrt{\phi}|\dx\dt \\
\leq &C\left(\int_{Q_{T}\backslash Q_{T,\eta}}\phi_k^{2\alpha}\dx\dt \right)^{1/2} + 
C\left(\int_{Q_{T}\backslash Q_{T,\eta}}\phi^{2\alpha}\dx\dt \right)^{1/2}
\end{align*}
and the strong convergence $\phi_k\to\phi$ in $L^2(Q_T)$ leads to
$$
\limsup_{k\to\infty}\mathcal{I}_{3,k}\leq  
C\left(\int_{Q_{T}\backslash Q_{T,\eta}}\phi^{2\alpha}\dx\dt \right)^{1/2}
\leq C\eta^{\alpha}.
$$
Putting the previous estimates on $\mathcal{I}_{1,k}$, $\mathcal{I}_{2,k}$, $\mathcal{I}_{3,k}$ together yields
$$
\limsup_{k\to\infty}\int_{Q_T}|\nabla\phi_k^{\alpha+1/2}-\nabla\phi^{\alpha+1/2}|\dx\dt\leq C\epsilon^{1-1/q} + 
C\eta^{\alpha}.
$$
The parameters $\epsilon>0$, $\eta>0$ are arbitrary. Hence, by choosing $\epsilon=\eta$ and letting $\eta\to 0$, we deduce that
$\nabla\phi_k^{\alpha+1/2}\to\nabla\phi^{\alpha+1/2}$ strongly in $L^1(Q_T)$.
Since $\nabla\phi_k^{\alpha+1/2}$ is bounded in $L^q(Q_T)$ with $q>2$,
it follows that $\nabla\phi_k^{\alpha+1/2}\to\nabla\phi^{\alpha+1/2}$ strongly in $L^2(Q_T)$. This finishes the proof of the lemma.
\end{proof}
     
    \section{Thermodynamic consistency based on a microforce balance} \label{app.thdyncons}
    To accurately describe the change in energy within an arbitrary control volume $V$, one has to consider energy changes due to macroscopic effects, i.e., work done by macroscopic stresses, as well as microscopic effects, i.e., changes due to the configuration and arrangement of atoms. For readability, we introduce the material derivative $\dot \psi = \partial_t \psi  + \vel \cdot \nabla \psi$. To describe the microscopic effects, we follow the approach of Gurtin~\cite{Gurtin1996, GuPoVi1996} and postulate a microforce balance law, independent of the macroscopic balance law. Specifically, Gurtin's microforce balance reads   
    \begin{align*}
        \int_{\partial V} \bm{\xi} \cdot \nu \dd s_x + \int_V (\pi + \gamma ) \dx = 0,  
    \end{align*}
    where $\bm{\xi} $ is a vector stress (microstress) and $\pi$ and $\gamma$ are scalar body forces, representing internal and external forces distributed over the volume. The atomic configuration is described by the phase field and order parameters $\phi_i, \ \Phi_v, \ i \in \mathcal{S}^*$. It is assumed that work due to microstresses is only performed with evolving fields, i.e., for an order parameter $\phi$, the working on the control volume due to microforces is given by $\int_{\partial V} ( \bm{\xi} \cdot \nu ) \dot{\phi} \dd s_x$. Let $V(t)$ denote an arbitrary volume in the bulk of the domain that is transported by $\vel$. We assume that we can split the microforces into $\bm{\xi} = \sum_{i \in \mathcal{S}^*}  \bm{\xi}_i + \bm{\xi}_v$ and define the working due to microscopic stresses as  
    \begin{align*}
         \sum_{i \in\mathcal{S}^*} \left( \int_{\partial V(t)} \dot{\phi}_i  \bm{\xi}_i \cdot \nu \dd s_x  \right) +  \int_{\partial V(t)} \dot{\Phi}_v  \bm{\xi}_v \cdot \nu \dd s_x.
    \end{align*}
    Considering the microscopic working, the total energy flux consists of the kinetic energy carried by diffusion, the working of the macroscopic stress, as well as the working of the microscopic stress, 
    \begin{align*}
     - \int_{\partial V(t)} \bm{j}_e \cdot \nu \dd s_x &= - \int_{\partial V(t)} \frac{|\vel|^2}{2} \bm{j} \cdot \nu \dd s_x +  \int_{\partial V(t)} ((\tilde \bT + \vel \otimes \bm{j}) \nu) \cdot \vel \dd s_x \\ & \quad + \sum_{i \in\mathcal{S}^*} \left( \int_{\partial V(t)} \dot{\phi}_i  \bm{\xi}_i \cdot \nu \dd s_x  \right) +  \int_{\partial V(t)} \dot{\Phi}_v  \bm{\xi}_v \cdot \nu \dd s_x .
    \end{align*}
    Following the Gibbs relation~\cite{GrootMazur62}, the entropy flux by diffusive mass transport is given by the transport of the chemical potentials $\mu_i$, $i \in \mathcal{S}^*$, i.e.,
    \begin{align*}
        \int_{\partial V(t)} \bm{j}_{\mathrm{ent}} \cdot \nu \dd s_x = - \sum_{i \in \mathcal{S}^*} \left( \int_{\partial V(t)} \frac{\mu_i}{\rhob_i} \bm{j}_i \cdot \nu \dd s_x \right) .
    \end{align*}
    Starting from Equation~\eqref{modeling.dissip.ineq}, one obtains the local balance equation 
    
    \begin{align*}
       & \partial_t f + \vel \cdot \nabla f  + \partial_t \left(\rho \frac{|\vel|^2}{2} \right) + \vel \cdot \nabla \left(\rho \frac{|\vel|^2}{2}\right) -\Div ((\tilde \bT^T + \bm{j} \otimes \vel ) \vel) \\ & \quad  - \sum_{i \in\mathcal{S}^*} \left( \Div (\dot{\phi}_i  \bm{\xi}_i) \right) - \Div (\dot{\Phi}_v  \bm{\xi}_v) + \sum_{i \in\mathcal{S}^*} \left( \Div \left( \frac{\mu_i}{\rhob_i} \bm{j}_i \right) \right) + \Div \left( \frac{|\vel|^2}{2} \bm{j} \right) = - \sigma 
    \end{align*}
    where $\sigma \geq 0$. Note that this formulation is equivalent to the dissipation inequality stated in~\cite{AbGaGr2012}. Recalling $f(\phi_{\mathcal{S}^*}, \Phi_v, \nabla \phi_{\mathcal{S}^*}, \nabla \Phi_v)$, we find the chain rule 
    \begin{align*}
        \dot{f} &= \partial_t f + \vel \cdot \nabla  f = \sum_{i \in \mathcal{S}^*} \left( \partial_{\phi_i} f \dot{\phi_i} + \left( \partial_{\nabla\phi_i} f \right) \cdot (\nabla \phi_i )^\cdot \right) + \left(\partial_{\Phi_v} f \right) \dot{\Phi}_v + \partial_{\nabla \Phi_v} \left( f \right) \cdot (\nabla \Phi_v )^\cdot. 
    \end{align*}
    Furthermore, expanding 
    \begin{align*}
        \nabla (\dot{\phi}) & = (\nabla \phi)^\cdot + (\nabla \vel )^T \nabla \phi, \\
        \Div (\dot{\phi}  \bm{\xi}) & = (\nabla \phi)^\cdot \cdot  \bm{\xi} + (\nabla \vel )  : (\nabla \phi \otimes   \bm{\xi}) + \dot{\phi}\Div  \bm{\xi} , 
    \end{align*}
    for any of the phase-field and order parameter, rewriting 
    \begin{align*}
         \partial_t\left( \rho \frac{|\vel|^2}{2}\right) + \vel \cdot \nabla \left( \rho \frac{|\vel|^2}{2}\right) &= \left( \partial_t \rho \right)\frac{|\vel|^2}{2} + \rho \vel \cdot \partial_t \vel + \vel \cdot (\nabla \rho) \frac{|\vel|^2}{2} + \vel \cdot \left( \rho (\nabla \vel ) \vel \right) \\
        &= \rho \left( \partial_t \vel + (\nabla \vel) \vel \right) \cdot \vel + \frac{|\vel|^2}{2}\left( \partial_t \rho + \vel \cdot \nabla \rho\right) \\
        &= \left( \Div \tilde \bT + \vel \Div \bm{j} \right)\cdot \vel - \frac{|\vel|^2}{2} \Div \bm{j} ,
    \end{align*}
    by the momentum balance equation~\eqref{eq_mombalance}, and substituting
    \begin{align*}
        \Div( \frac{\mu_i}{\rhob_i} \bm{j}_i) & = \frac{\mu_i}{\rhob_i} \Div \bm{j}_i + \nabla  \frac{\mu_i}{\rhob_i} \cdot \bm{j}_i =  -\mu_i \dot{\phi_i} + \nabla \frac{\mu_i}{\rhob_i}  \cdot \bm{j}_i  
    \end{align*}
    using the volume balance equations~\eqref{eq_volumebalance}, simplifies the equation to 
        \begin{align*}
        - \sigma & = \left(\sum_{i \in\mathcal{S}^*} \left( \partial_{\phi_i} f - \mu_i  - \Div  \bm{\xi}_i \right) \dot{\phi_i} \right) + \left( \sum_{i \in\mathcal{S}^*} \left( \partial_{\nabla \phi_i} f -  \bm{\xi}_i \right) \cdot  (\nabla \phi_i)^\cdot  \right) \\
        & \quad +  \left(\partial_{\Phi_v} f - \Div \bm{\xi}_v\right) \dot{\Phi}_v +\left( \partial_{\nabla \Phi_v} f - \bm{\xi}_v \right) \cdot  (\nabla \Phi_v)^\cdot \\
        & \quad - (\nabla \vel)^T : \left( \sum_{i \in\mathcal{S}^*} \nabla \phi_i \otimes  \bm{\xi}_i + \nabla \Phi_v \otimes  \bm{\xi}_v \right) + \Div \tilde \bT \cdot \vel - \Div(\tilde \bT^T \vel)   \\
        & \quad + \vel \Div \bm{j} \cdot \vel - \Div ((\vel \otimes \bm{j} )^T \vel) - \frac{|\vel|^2}{2} \Div \bm{j} + \Div (\frac{|\vel|^2}{2}\bm{j} ) + \sum_{i \in\mathcal{S}^*} \nabla \frac{\mu_i}{\rhob_i} \cdot \bm{j}_i .
    \end{align*}
    Finally, rewriting 
    \begin{align*}
        \Div \tilde \bT \cdot \vel - \Div (\tilde \bT^T \vel) &= - \tilde \bT : (\nabla \vel), \\
        \vel \Div \bm{j} \cdot \vel - \Div ((\vel \otimes \bm{j})^T \vel) &= - ((\nabla \vel) \bm{j}) \cdot \vel - (\vel \otimes \bm{j}) : (\nabla \vel), \\
        \Div \left(\frac{|\vel|^2}{2} \bm{j} \right) &= \frac{|\vel|^2}{2}\Div \bm{j} + ((\nabla \vel) \bm{j}) \cdot \vel ,
        \end{align*}
    and recalling $\bS = \tilde \bT + \vel \otimes \bm{j} + p \bI$ and $ p \bI : (\nabla \vel) = p \left(  \Div \vel  \right)= 0$ yields
    \begin{align*} 
        - \sigma & = \left(\sum_{i \in\mathcal{S}^*} \left( \partial_{\phi_i} f - \mu_i - \Div \bm{\xi}_i \right) \dot{\phi_i} \right) +  \left(\partial_{\Phi_v} f - \Div \bm{\xi}_v\right) \dot{\Phi}_v \\ 
        & \quad + \left( \sum_{i \in\mathcal{S}^*} \left( \partial_{\nabla \phi_i} f - \bm{\xi}_i \right) \cdot  (\nabla \phi_i)^\cdot  \right) +  \left( \partial_{\nabla \Phi_v} f - \bm{\xi}_v \right) \cdot  (\nabla \Phi_v)^\cdot\\
        & \quad - \left( \bS +  \left( \sum_{i \in\mathcal{S}^*} \nabla \phi_i \otimes \bm{\xi}_i \right) + \nabla \Phi_v \otimes \bm{\xi}_v \right) : D ( \vel)   
        \\ & \quad - \left( \sum_{i \in\mathcal{S}^*} \nabla \phi_i \otimes \bm{\xi}_i + \nabla \Phi_v \otimes \bm{\xi}_v\right) : W(\vel)^T  + \left( \sum_{i \in\mathcal{S}^*} \nabla \frac{\mu_i}{\rhob_i} \cdot \bm{j}_i \right). 
    \end{align*}
    At this point we select constitutive equations compatible with the condition $\sigma \geq 0$. Under the assumption that $\dot{\phi}_i$, $ (\nabla \phi_i) ^\cdot$, $\dot{\Phi}_v$, $\left( \nabla \Phi_v \right)^\cdot$, $D(\vel )$, $W( \vel ) $, and $\mu_i$ can take arbitrary values, we require all terms to obey the nonnegativity constraint separately. Starting with the third and fourth terms, notice that $\partial_{\nabla \phi_i} f$ is independent of $ (\nabla \phi_i)^\cdot$. Therefore, setting 
    \begin{align*}
         \partial_{\nabla \phi_i} f =  \bm{\xi}_i \text{ for } i \in \mathcal{S}^* \text{ and } 
         \partial_{\nabla \Phi_v} f = \bm{\xi}_v
    \end{align*}
    ensures these terms vanish. This implies the constitutive equations $ \bm{\xi}_i(\phi_{\mathcal{S}^*}, \Phi_v, \nabla \phi_{\mathcal{S}^*}, \nabla \Phi_v)$ and $ \bm{\xi}_v(\phi_{\mathcal{S}^*}, \Phi_v, \nabla \phi_{\mathcal{S}^*}, \nabla \Phi_v)$. Furthermore, it is sufficient to require $\bm{\xi}_i \otimes \nabla \phi_i = \nabla \phi_i \otimes  \bm{\xi}_i$ for the second-to-last term to vanish. Equating $\left( \bm{\xi}_i \otimes \nabla \phi_i \right) \nabla \phi_i = \left( \nabla \phi_i \otimes  \bm{\xi}_i \right) \nabla \phi_i$, we can conclude $(\bm{\xi}_i \cdot \nabla \phi_i ) \nabla \phi_i = |\nabla \phi_i|^2  \bm{\xi}_i$. It follows that $ \bm{\xi}_i$ is parallel to $\nabla \phi_i$, and we can write  
    \begin{align*}
         \bm{\xi}_i = \alpha_i \nabla \phi_i = \partial_{\nabla \phi_i} f
    \end{align*}
    for some (phase-field-dependent) scalar function $\alpha_i$. Similarly, $\bm{\xi}_v = \alpha_v \nabla \Phi_v = \partial_{\nabla \Phi_v} f$. Consequently, to cancel out the first term, the chemical potentials take the form %
    \begin{align*}
       \mu_i =  \partial_{\phi_i} f - \Div \partial_{\nabla \phi_i} f . 
    \end{align*}
    With these restrictions, the dissipation inequality simplifies to 
    \begin{align*}
        - \sigma  & =  \left(\partial_{\Phi_v} f - \Div \partial_{\nabla \Phi_v} f\right) \dot{\Phi}_v   - \left( \bS +  \left( \sum_{i \in\mathcal{S}^*} \nabla \phi_i \otimes  \bm{\xi}_i \right) + \nabla \Phi_v \otimes  \bm{\xi}_v \right) : D (\vel)  \\ & \quad + \left( \sum_{i \in\mathcal{S}^*} \nabla \frac{\mu_i}{\rhob_i} \cdot \bm{j}_i \right) \leq 0 .
    \end{align*}
    It remains to specify the other constitutive equations to guarantee thermodynamic consistency. We state linear phenomenological laws, ensuring a nonnegative contribution for each term individually. First, we employ Newton's rheological law and find the relation $\bS + \sum_{i \in\mathcal{S}^*} \nabla \phi_i \otimes \partial_{\nabla \phi_i} f+ \nabla \Phi_v \otimes \partial_{\nabla \Phi_v} f = 2 \eta(\phi_{\mathcal{S}^*})D (\vel)$, with $\eta(\phi_{\mathcal{S}^*}) \geq 0$ the (phase-dependent) viscosity parameter. 
    To ensure nonnegativity of the last term, we recall the constraint $\sum_i \bm{J}_i = 0$ and note that the mass fluxes cannot be chosen independently. Following Fick's diffusion law, we set $\bm{j}_i = - \sum_j m_{ij} \nabla \mu_j$. To ensure nonnegativity, we require the matrix with entries $(\bM)_{ij} = \frac{m_{ij}}{\rhob_i}$ to be symmetric positive semidefinite. Furthermore, we require $\sum_i \frac{m_{ij}}{\rhob_i} = 0$ to ensure the volume constraint. Lastly, we define the evolution of the order parameter   
    \begin{align*}
        \dot{\Phi}_v = - m_v(\phi_{\mathcal{S}^*}) \left(\partial_{\Phi_v} f - \Div \partial_{\nabla \Phi_v} f \right ) 
    \end{align*}
    for a scalar $m_v(\phi_{\mathcal{S}^*}) \geq 0$, which may depend on the volume fractions $\phi_i$. These assumptions lead to the entropy production 
    \begin{align*}
        \sigma  = m_v(\phi_{\mathcal{S}^*}) |\partial_{\Phi_v} f - \Div \partial_{\nabla \Phi_v} f|^2 + 2 \eta(\vec{\phi}) |D (\vel)|^2 + \sum_{i, j \in \mathcal{S}^*} \nabla \mu_i \cdot  \frac{m_{ij}(\phi_{\mathcal{S}^*}, \Phi_v)}{\rhob_i} \nabla \mu_j 
    \end{align*}
    and the resulting energy balance reads
    \begin{align*}
        & \frac{\dd}{\dt} \int_{V(t)} \rho \frac{|\vel|^2}{2} + f \dx - \int_{\partial V(t)} (\tilde \bT + \vel \otimes \bm{j}) \nu  \cdot \vel \dd s_x - \sum_{i \in\mathcal{S}^*} \left( \int_{\partial V(t)} \dot{\phi_i}  \bm{\xi}_i \cdot \nu \dd s_x \right) \\ 
        & \quad -\int_{\partial V(t)} \dot{\Phi}_v  \bm{\xi}_v\cdot\nu \dd s_x 
        + \sum_{i \in\mathcal{S}^*} \left( \int_{\partial V(t)} \frac{\mu_i}{\rhob_i} \bm{j}_i \cdot \nu \dd s_x \right) + \int_{\partial V(t)} \frac{|\vel|^2}{2} \bm{j} \cdot \nu \dd s_x  \\ 
        & \quad + \int_{V(t)} \left( m_v(\phi_{\mathcal{S}^*}) |\partial_{\Phi_v} f - \Div \partial_{\nabla \Phi_v} f|^2 + 2 \eta(\vec{\phi}) |D (\vel)|^2 + \sum_{i, j \in \mathcal{S}^*} \nabla \mu_i \frac{m_{ij}(\phi_{\mathcal{S}^*}, \Phi_v)}{\rhob_i} \nabla \mu_j  \right) \dx 
        \\ & \quad = 0 . 
    \end{align*}
	
	\section{Derivation of the nondimensional equations} \label{app.nondimeq}
    We fix the characteristic length, time, and mass scales $x_c$, $t_c$, and $m_c$. Then $\nabla = \frac{1}{x_c} \overset{*}{\nabla}$, $\frac{\partial}{\partial t} = \frac{1}{t_c} \frac{\partial}{\partial t^*}$ and 
    \begin{align*}
    \bx &= \bx^* x_c , &  t&= t^* t_c , \  &\phi_i & = \phi_i^*, \  & \Phi_v & = \Phi_v^* ,\ & \rho  &= \rho^* \frac{m_c}{x_c^3} ,\ & \vel &= \vel^* \frac{x_c}{t_c}, \\
    \mu_i &= \mu_i^* \frac{m_c}{x_c t_c^2} ,\ & f &= f^*  \frac{m_c}{x_c t_c^2} ,\ & p &= p^*\frac{m_c}{x_c t_c^2}, \ & \bm{j} &= \bm{j}^*\frac{m_c}{x_c^2 t_c}, \ & \bm{J} &= \bm{J}^*\frac{x_c}{t_c},
    \end{align*}
    with the dimensionless $(\cdot)^*$ and characteristic $(\cdot)_c$. In particular, setting $m_c = \frac{R \theta}{V^0}x_c t_c^2$ simplifies the nondimensional local free energy density (compare with equation~\eqref{eq_freeenergy}) to
    \begin{align*}
        f^{\mathrm{loc}*} &= \left( 1- \omega(\Phi_v)\right) \left( \sum_{i \in \mathcal{S}^*} \frac{\phi_i}{N_i} \log \phi_i + \sum_{i, j \in \mathcal{S}^*} \chi_{ij}^* \phi_i \phi_j  \right) +  \omega(\Phi_v) \left( \sum_{i \in \mathcal{S}^*} \frac{\phi_i}{N_i} \log \frac{\phi_i}{\phi^{\mathrm{sat}}_i}\right),
    \end{align*}
    where $\chi_{ij}^* = \frac{\chi_{ij}}{2}$, where we symmetrically extended $\chi_{ij}$ and set $\chi_{ii} = 0$. Combining the momentum balance~\eqref{eq_momentumbalance} with the constitutive equation~\eqref{eq_const_S} and the definition of the nonlocal free energy density yields, after rescaling 
    \begin{align*}
        &\rho^*\pa_{t^*}\vel^* + \left(\rho^*\vel^* 
		+ \sum_{i \in \mathcal{S}^*} \bm{j}_i^*\right)\cdot\overset{*}{\nabla}\vel^*
		& \\
		& =\overset{*}{\Div}\left(2 \frac{x_c t_c }{m_c} \eta D^*(\vel^*)\right)
		-\sum_{i \in \mathcal{S}^*} \frac{t_c^2}{x_c m_c} \epsilon_i^2 \overset{*}{\Div}( \overset{*}{\nabla} \phi_i^*\otimes  \overset{*}{\nabla} \phi_i^*)
		- \frac{t_c^2}{x_c m_c} \lambda^2 \overset{*}{\Div}(\overset{*}{\nabla}\Phi_v^* \otimes\overset{*}{\nabla}\Phi_v^*) - \overset{*}{\nabla} p^* .
    \end{align*}
    Similarly, it follows from Equation~\eqref{eq_const_j} and the relation $\bm{J}_i = \frac{1}{\rhob_i} \bm{j}_i$ that
    \begin{align*}
        \bm{J}^*_i =  -\sum_{j \in \mathcal{S}^*} \frac{m_c}{x_c^3 t_c} \frac{m_{ij}}{\rhob_i} \overset{*}{\nabla} \mu_j^*.  
    \end{align*}
    From Equation~\eqref{eq_const_mu} and the free energy density follows the nondimensional chemical potential, 
    \begin{align*}
        \mu_i^* = \pa_{\phi_i^*} f^{\mathrm{loc}*} - \frac{t_c^2 }{x_c m_c}  \epsilon_i^2 \overset{*}{\Delta} \phi_i^* .
    \end{align*}
    Correspondingly, following Equation~\eqref{eq_const_Phiv}, 
    \begin{align*}
        \pa_{t^*} \Phi_v^* + \vel^* \cdot \overset{*}{\nabla} \Phi_v ^* = - \frac{m_c}{x_c t_c} m_v \left (\pa_{\Phi_v^*}f^{\mathrm{loc}*} - \frac{t_c^2}{x_c m_c} \lambda^2 \overset{*}{\Delta} \Phi_v^* \right) .
    \end{align*}
    The scaling of the remaining governing equations is straightforward, and their derivations are omitted here. We introduce the nondimensional viscosity and densities 
    \begin{align*}
        \eta^* = \frac{x_c t_c}{m_c}\eta \text{ and }  \rhob_i^* = \frac{x_c^3}{m_c} \rhob_i, 
    \end{align*}
    and furthermore introduce the parameters
    \begin{align*}
        \beta_i & = \frac{t_c^2}{x_cm_c}\epsilon^2_i, \ &\beta_{v} &= \frac{t_c^2}{x_c m_c} \lambda^2, \ &\gamma_{ij} & = \frac{m_c}{x_c^3t_c \rhob_i} m_{ij} , \ &\gamma_{v} & =  \frac{m_c}{x_c t_c} m_{v}.  
    \end{align*}
  Substituting these parameters and dropping the superscripts yields the nondimensional governing equations.

\bibliographystyle{plain}
\bibliography{references}

\end{document}